\documentclass[11pt,a4paper]{article}
\usepackage[utf8]{inputenc}
\usepackage{amsmath,amssymb,amsthm,booktabs,array}
\usepackage[margin=2.3cm]{geometry}
\usepackage{tikz}
\usepackage{hyperref}
\theoremstyle{plain}
\newtheorem{theorem}{Theorem}[section]
\newtheorem{lemma}[theorem]{Lemma}
\newtheorem{proposition}[theorem]{Proposition}
\newtheorem{corollary}[theorem]{Corollary}
\theoremstyle{definition}
\newtheorem{definition}[theorem]{Definition}
\newtheorem{example}[theorem]{Example}
\newtheorem{remark}[theorem]{Remark}
\newcommand{\N}{\mathbb{N}}
\newcommand{\PG}{\mathrm{PG}}
\newcommand{\PM}{\mathrm{PM}}
\newcommand{\PF}{\mathrm{PF}}
\newcommand{\cs}{C}
\newcommand{\SSk}{\Sigma(\kappa)+\Sigma(\kappa)}

\title{\textbf{Creation Sequences for Numerical Semigroups\\ and Exact Counts of Semigroup Partitions}}
\author{Rameez Raja\thanks{Department of Mathematics, National Institute of Technology Srinagar, Srinagar 190006,
Jammu and Kashmir, India. Emails: \texttt{rameeznaqash@nitsri.ac.in} (R.~Raja), \texttt{zeenat\string_2024phamth005@nitsri.ac.in} (Z.~Rashid).} \and Zeenat Rashid\footnotemark[1]}
\date{}

\begin{document}
\maketitle

\begin{abstract}
A numerical semigroup is a subset of the non-negative integers that contains $0$, is closed under addition and has
finite complement. A threshold graph is built by adding vertices one at a time, each new vertex being isolated or
dominating, and the record of these choices is its \emph{creation sequence}. Reading the creation sequence as a
lattice path gives an integer partition and, through the Keith--Nath correspondence, a numerical set in which the
dominating vertices are the elements and the isolated vertices are the gaps. We show that threshold graphs with an
isolated vertex correspond bijectively to numerical sets, and we study numerical semigroups through their creation
sequences. We prove that in the creation sequence of a numerical semigroup of multiplicity $m$ every block of ones
and every non-initial block of zeros has length less than $m$, we give criteria for partitions with a unique smallest
part and with two distinct part sizes, and we classify all partitions with at most four parts. We show that every
partition $\lambda$ with $g$ parts corresponding to a numerical semigroup satisfies $\ell(\lambda)\ge\lambda_1$ and,
more strongly, fits inside the staircase $(g,g-1,\dots,1)$; equivalently, its creation sequence begins with a Dyck
path. Its size equals the genus plus the Weierstrass weight of the semigroup, and asymptotically at most half of all
partitions correspond to numerical semigroups. Finally, we obtain exact formulas for the numbers of partitions of $n$
corresponding to numerical semigroups of genus $g$, multiplicity $m$ and Frobenius number $f$, valid in the sharp
ranges $g\ge n/2$, $m>n/3$ and $f\ge n-\frac{n-1}{4}$. These formulas solve the open problems posed by Burson, Nam and
Sisneros-Thiry \cite{BNS}, recover all of their counting theorems, and correct several of their statements.
\end{abstract}



%======================================================================
\section{Introduction}
%======================================================================

\subsection{Numerical semigroups}

A \emph{numerical semigroup} is a subset $S$ of the non-negative integers $\N_0$ that contains $0$, is closed under
addition, and has finite complement in $\N_0$. Numerical semigroups arise naturally in the study of linear Diophantine
equations: for coprime positive integers $a_1,\dots,a_k$, the set of integers $d$ for which $a_1x_1+\dots+a_kx_k=d$ has
a solution in non-negative integers is the numerical semigroup $\langle a_1,\dots,a_k\rangle$ generated by
$a_1,\dots,a_k$. The largest integer that is not representable in this way is the \emph{Frobenius number}. The number
of non-representable integers is the \emph{genus}, and the smallest positive representable integer is the
\emph{multiplicity}. Sylvester showed that $\langle a,b\rangle$ has Frobenius number $ab-a-b$ and genus
$(a-1)(b-1)/2$. No such formula is known for three or more generators, and the \emph{Frobenius problem} of computing
these invariants has a long history \cite{RGS,AGS}. Numerical semigroups also appear as value semigroups of curve
singularities and as Weierstrass semigroups of points on algebraic curves, where the genus of the semigroup is the
genus of the curve and the weight of the semigroup measures how special the point is \cite{EH}. Among the many
distinguished classes of numerical semigroups, the symmetric, pseudo-symmetric, almost symmetric and Arf semigroups
play a central role in commutative algebra, since they correspond to Gorenstein, Kunz, almost Gorenstein and Arf
semigroup rings.

A fruitful combinatorial point of view encodes a numerical semigroup by an integer partition. Keith and Nath
\cite{KN}, motivated by the study of partitions with prescribed hook lengths and by the hook length formula for the
dimensions of irreducible representations of symmetric groups, observed the following. If one walks along
$0,1,2,\dots$ and takes an east step for each element of a numerical set $S$ and a north step for each gap, the
resulting lattice path is the profile of a Young diagram. In this way every numerical set determines a partition
$\lambda(S)$, and the restriction of this map to numerical semigroups is injective. Under this correspondence: the genus of $S$ is the number of parts of $\lambda(S)$; the Frobenius number of $S$ is the largest hook length of $\lambda(S)$; the multiplicity of $S$ is one more than the number of parts of $\lambda(S)$ equal to $1$; the gaps of $S$ are exactly the hook lengths in the first column of $\lambda(S)$.

Constantin, Houston-Edwards and Kaplan \cite{CHK} used this correspondence to relate numerical sets to core partitions
and to integer points in polytopes. They showed that the partitions with largest hook $N$ that come from numerical
semigroups form a vanishing proportion of all partitions with largest hook $N$, and they conjectured that the same
holds when partitions are counted by size.

A series of papers by Tutaş and her coauthors develops the partition
side of the theory for Arf semigroups. Tutaş, Karakaş and Gümüşbaş \cite{TKG19} introduced \emph{Arf partitions} and
characterised Arf numerical semigroups through the hook sets of the Young tableaux of the associated partitions.
Tutaş \cite{Tut19} extended this work, studying the relationship between Arf partitions, Young dual diagrams and
rational Young diagrams. She showed that the partitions arising from gap sets of numerical semigroups form a semigroup
under an intersection operation, and she related Arf partitions to semigroup polynomials, obtaining generating
functions for families such as the Lipman semigroups. Karakaş and Tutaş \cite{KT20} gave a decomposition of partitions
and of numerical sets that is compatible with the Keith--Nath correspondence. The authors in 
\cite{GTE20} introduced \emph{almost symmetric Arf partitions} and used partition theory to count almost symmetric Arf
semigroups of a given genus.

Süer and Yeşil \cite{SY21} studied the Young
diagrams of symmetric and pseudo-symmetric numerical semigroups. They introduced new operations on Young diagrams and
on numerical semigroups that yield decompositions of these objects, and they determined when the associated semigroup
rings have Gorenstein and Kunz subrings. Süer and Yeşil \cite{SY24} investigated special subdiagrams of Young
diagrams and their numerical semigroups. More recently, Yeşil \cite{Yes25} proved that every almost symmetric numerical
semigroup decomposes uniquely, at the level of Young diagrams, into a numerical semigroup, its dual, and a further
numerical semigroup. In these works the dual (conjugate) of a Young diagram plays a central role; in particular,
symmetric numerical semigroups correspond to self-conjugate diagrams. We use the same phenomenon in
Section~\ref{sec:sym}. Related partition-theoretic viewpoints appear in recent work of Campioni
\cite{Cam26} on unrefinable partitions into distinct parts and numerical semigroups, and of Aragona, Civino and
Campioni \cite{ACC26}, who characterise maximal unrefinable partitions geometrically via the Keith--Nath
transformation and Young diagrams. Burson, Nam and Sisneros-Thiry \cite{BNS} refined the question of
Constantin, Houston-Edwards and Kaplan by counting
$
\PG(n,g), \PF(n,f), ~\text{and}~ \PM(n,m),
$
the numbers of partitions of $n$ corresponding to numerical semigroups of genus $g$, Frobenius number $f$ and
multiplicity $m$. They proved exact identities along lines close to the extreme cases, such as
$\PG(n,g)=p(n-g)$ for $g\ge\frac23n$ and $\PF(n,n)=\lceil n/2\rceil$. They left three problems open: finding the
functions that count the exceptional partitions for the genus (\cite[Remark 1]{BNS}) and for the multiplicity
(\cite[Remark 4]{BNS}), and finding uniform theorems valid for every shift (\cite[Section 6]{BNS}).

A \emph{threshold graph} is a graph that can be built from a single vertex by repeatedly adding either an isolated
vertex or a dominating vertex \cite{HIS,MP}. Threshold graphs have many equivalent characterisations, for instance by
vertex weights and a threshold, or by forbidden induced subgraphs. They are determined by their degree sequences. The
construction of a threshold graph is recorded by a binary word, its \emph{creation sequence}, with $0$ for an isolated
vertex and $1$ for a dominating vertex. Up to isomorphism, threshold graphs on $N$ vertices correspond bijectively to
binary words of length $N-1$.

In this paper we connect threshold graphs with numerical semigroups through their creation sequences. We associate
numerical semigroups with threshold graphs in three steps:
\[
\text{threshold graph } G\ \longmapsto\ \text{creation sequence } s(G)\ \longmapsto\ \text{partition }\lambda_G
\ \longmapsto\ \text{numerical set } S_G .
\]
The partition $\lambda_G$ is obtained by reading $s(G)$ as a lattice path ($1$ = east, $0$ = north), and $S_G$ is the
Keith--Nath numerical set of $\lambda_G$. Thus the vertex added at step $i$ is dominating exactly when $i\in S_G$.
Threshold graphs with an isolated vertex correspond bijectively to numerical sets (Theorem~\ref{thm:graphbij}), and
graph-theoretic quantities translate into semigroup invariants (Proposition~\ref{prop:graphdict}). The creation
sequence makes additive closure visible: $S_G$ fails to be a semigroup exactly when the vertices added at two steps
$x$ and $y$ are dominating while the vertex added at step $x+y$ is isolated.

Section~\ref{sec:prelim} collects preliminaries on numerical semigroups, partitions and the Keith--Nath
correspondence. Section~\ref{sec:cs} introduces threshold graphs and creation sequences and proves the correspondence
theorem. Sections~\ref{sec:struct}--\ref{sec:tree} contain the structural results, and Section~\ref{sec:count} the
counting results. Section~\ref{sec:solutions} collects the solutions of the open problems of \cite{BNS}, and the last
section lists further questions. More precisely, working entirely with creation sequences, we obtain the following interesting results.
\begin{enumerate}
\item \emph{Structural criteria} (Section~\ref{sec:struct}). In the creation sequence of a numerical semigroup of
multiplicity $m$, every block of ones and every non-initial block of zeros has length less than $m$
(Theorem~\ref{thm:blocks}). In partition language, no part occurs more often than the smallest part, and consecutive
distinct parts differ by at most that multiplicity. We give complete criteria for partitions with a unique smallest part
(Theorem~\ref{thm:unique}) and with two distinct part sizes (Theorem~\ref{thm:twosizes}), and a column criterion
deciding every creation sequence (Theorem~\ref{thm:column}).
\item \emph{Classification by patterns} (Section~\ref{sec:patterns}). We classify all partitions with at most four
parts (Table~\ref{tab:small}).
\item \emph{Symmetry and density} (Section~\ref{sec:sym}). Conjugation of partitions reverses and complements the
creation sequence. Consequently every semigroup partition satisfies $\ell(\lambda)\ge\lambda_1$, with equality
exactly for self-conjugate partitions, that is, symmetric semigroups (Theorem~\ref{thm:ellbound}). We deduce that
asymptotically at most half of all partitions correspond to numerical semigroups (Theorem~\ref{thm:density}).
\item \emph{The staircase theorem} (Section~\ref{sec:staircase}). Every semigroup partition with $g$ parts lies inside
the staircase $(g,g-1,\dots,1)$; equivalently, its creation sequence satisfies a ballot condition and begins with a
Dyck path (Theorem~\ref{thm:staircase}, Corollary~\ref{cor:catalan}). The size of the partition equals the genus plus
the weight of the semigroup (Corollary~\ref{cor:weight}).
\item \emph{The semigroup tree} (Section~\ref{sec:tree}). Deleting the largest part of a semigroup partition gives a
semigroup partition; this realises the tree of numerical semigroups in partition language.
\item \emph{Exact counting and the open problems of \cite{BNS}} (Sections~\ref{sec:count} and \ref{sec:solutions}).
We prove
\begin{align*}
\PG(n,g)&=p(n-g)-E(2n-3g) && (2g\ge n),\\
\PM(n,m)&=\bar p(n-m+1)-H(n-2m) && (m-1\le n<3m),\\
\PF(n,n-k)&=p(k)\lceil\tfrac{n-k}{2}\rceil-C_k(n-k) && (n\ge4k+1),
\end{align*}
where $E$, $H$ and $C_k$ are explicit finite counts of tail partitions defined through their creation sequences. We
show that the three ranges are sharp, and we give exact formulas on the boundary lines. These results solve the open
problems of \cite[Remarks 1 and 4, Section 6]{BNS} and recover all of their counting theorems; in the process we
correct their Theorem 3, Corollary 4, Theorem 7(3)-(5) and Theorem 10-(4).
\end{enumerate}

%======================================================================
\section{Preliminaries}\label{sec:prelim}
%======================================================================

\subsection{Numerical sets and numerical semigroups}

Let $\N_0=\{0,1,2,\dots\}$. A \emph{numerical set} is a subset $S\subseteq\N_0$ with $0\in S$ and $\N_0\setminus S$
finite. A numerical set is a \emph{numerical semigroup} if it is closed under addition. The elements of
$\N_0\setminus S$ are the \emph{gaps} of $S$, and we record the following invariants:

(i) the \emph{genus} $g(S)=|\N_0\setminus S|$;

(ii) the \emph{Frobenius number} $F(S)=\max(\N_0\setminus S)$, with the convention $F(\N_0)=-1$;

(iii) the \emph{conductor} $c(S)=F(S)+1$, the least integer $c$ with $[c,\infty)\subseteq S$;

(iv) the \emph{multiplicity} $m(S)=\min(S\setminus\{0\})$.

For $A\subseteq\N_0$ we write $\langle A\rangle$ for the set of all finite sums of elements of $A$; if $\gcd(A)=1$ this
is a numerical semigroup. Every numerical semigroup $S$ has a unique minimal system of generators, consisting of its
\emph{minimal generators}: the nonzero elements that are not the sum of two nonzero elements of $S$. Their number is
the \emph{embedding dimension}, and it is at most $m(S)$. For $0\ne n\in S$, the \emph{Apéry set}
$\mathrm{Ap}(S,n)=\{s\in S:\ s-n\notin S\}$ contains exactly one element in each residue class modulo $n$.

\begin{example}\label{ex:basic}
The numerical semigroup $\langle3,5\rangle=\{0,3,5,6,8,9,10,\dots\}$ has gaps $\{1,2,4,7\}$, so $g=4$, $F=7$ and $m=3$. Its Apéry set with
respect to $3$ is $\{0,5,10\}$. In accordance with Sylvester's formulas, $F=3\cdot5-3-5$ and $g=(3-1)(5-1)/2$. For $\langle3,5,7\rangle=\N_0\setminus\{1,2,4\}=\{0,3,5,6,7,\dots\}$, $g=3$, $F=4$ and $m=3$.
Note that the numerical set $\N_0\setminus\{1,3,4\}=\{0,2,5,6,\dots\}$ is not a semigroup, since $2+2=4$ is a gap.

\end{example}

A numerical semigroup $S$ is \emph{symmetric} if, for all $x\in\mathbb Z$,
$x\in S$ if and only if $F(S)-x\notin S$. It is \emph{pseudo-symmetric} if $F(S)$ is even and, for all
$x\ne F(S)/2$, $x\in S$ if and only if $F(S)-x\notin S$. It is \emph{irreducible} if it cannot be written as the
intersection of two numerical semigroups properly containing it. A numerical semigroup is irreducible if and only if it
is symmetric or pseudo-symmetric, and this happens if and only if $g(S)=\lceil(F(S)+1)/2\rceil$ \cite{RGS}. For every
numerical semigroup, $g(S)\ge(F(S)+1)/2$. A numerical semigroup is \emph{Arf} if $x+y-z\in S$ for all $x\ge y\ge z$ in
$S$ \cite{RGS}; Arf semigroups and their partitions are the subject of \cite{TKG19,Tut19,GTE20}.

\begin{example}
$\langle3,5\rangle$ and $\langle4,5,6\rangle$ are symmetric: both have $g=4=(F+1)/2$ with $F=7$. For
$\langle4,5,6\rangle$, for instance, the gaps $1,2,3,7$ are exactly the numbers $7-s$ with $s\in\{6,5,4,0\}$.
$\langle3,5,7\rangle$ is pseudo-symmetric: $F=4$ is even and $g=3=(F+2)/2$. $\langle4,5,6,7\rangle=\{0,4,5,\dots\}$ has
$F=3$ and $g=3>\lceil(3+1)/2\rceil$, so it is not irreducible.
\end{example}

If $S\ne\N_0$ is a numerical semigroup, then $S\cup\{F(S)\}$ is again a numerical
semigroup, of genus $g(S)-1$. Therefore all numerical semigroups form a rooted tree with root $\N_0$, in which the parent of
$S$ is $S\cup\{F(S)\}$. The children of $S$ are the semigroups $S\setminus\{x\}$ for the minimal generators $x>F(S)$
\cite{RGS}. The number $n_g$ of numerical semigroups of genus $g$ satisfies
$2F_g\le n_g\le1+3\cdot2^{g-3}$ for $g\ge3$, with $F_g$ the Fibonacci numbers \cite{BA}, and $n_g$ grows like $C\varphi^g$,
where $\varphi$ is the golden ratio \cite{Zhai}.

\subsection{Partitions and Young diagrams}

A \emph{partition} of $n$ is a non-increasing sequence $\lambda=(\lambda_1,\dots,\lambda_k)$ of positive integers
with $|\lambda|=\sum\limits_{i=1}^k\lambda_i=n$. We write $\lambda\vdash n$, $\ell(\lambda)=k$ for the number of parts, and use
exponents for repeated parts, that is, $(3,2^2,1^3)=(3,2,2,1,1,1)$. The \emph{Young diagram} of $\lambda$ is the
left-justified array of cells with $\lambda_i$ cells in row $i$, rows numbered from the top. The \emph{conjugate}
$\lambda'$ is obtained by reflecting the diagram across its main diagonal, so that
$\lambda'_j=\#\{i:\lambda_i\ge j\}$; a partition is \emph{self-conjugate} if $\lambda=\lambda'$. The \emph{hook} of a
cell $(i,j)$ consists of the cell itself, the cells to its right in row $i$ (the \emph{arm}), and the cells below it
in column $j$ (the \emph{leg}). Its \emph{hook length} is
\[
h_{(i,j)}(\lambda)=\lambda_i-j+\lambda'_j-i+1 .
\]
In particular $h_{(1,1)}(\lambda)=\lambda_1+\ell(\lambda)-1$. A partition is a hook if it has the form
$(a,1^b)$. For partitions $\mu,\lambda$ we write $\mu\subseteq\lambda$ if $\mu_i\le\lambda_i$ for all $i$. We use the following counting functions:
\[
p(n)=\#\{\lambda\vdash n\},\qquad \bar p(n)=\#\{\lambda\vdash n:\ 1\notin\lambda\},\qquad
a(n)=\bar p(n)-1\ (n\ge2),
\]
\[
\bar a(n)=\#\{\lambda\vdash n:\ 1\notin\lambda,\ \ell(\lambda)\ge2,\ \lambda_1=\lambda_2\},\qquad
b(n)=\#\{\lambda\vdash n+2:\ \ell(\lambda)\le n,\ \lambda_1\le n\}.
\]
Removing one part equal to $1$ shows that $p(n)-p(n-1)=\bar p(n)$. By the Hardy--Ramanujan formula \cite{HR},
$p(n)\sim\frac{1}{4n\sqrt3}e^{\pi\sqrt{2n/3}}$; in particular $p(n-1)/p(n)\to1$.

\subsection{The Keith--Nath correspondence}

Let $S\ne\N_0$ be a numerical set with Frobenius number $F$. Walk along $0,1,\dots,F$, taking an east step for each
element of $S$ and a north step for each gap. Starting from the bottom-left corner, this traces the profile of a Young
diagram, and we denote the corresponding partition by $\lambda(S)$. The rows of $\lambda(S)$ correspond to the north
steps, that is, to the gaps. The length of the row belonging to a gap $x$ is the number of elements of $S$ smaller than
$x$. Conversely, every partition $\lambda\ne\emptyset$ equals $\lambda(S)$ for exactly one numerical set $S$
\cite{KN}, and we set $\lambda(\N_0)=\emptyset$. We write $S^{KN}(\lambda)$ for the inverse map.

\begin{proposition}[\cite{KN,BNS}]\label{prop:KNdict}
Let $S\ne\N_0$ be a numerical set and $\lambda=\lambda(S)$. Then $g(S)=\ell(\lambda)$ and $F(S)=h_{(1,1)}(\lambda)$.
The gaps of $S$ are exactly the hook lengths $h_{(i,1)}(\lambda)$ of the first column. If $S$ is a numerical semigroup
and $S\ne\N_0$, then $m(S)-1$ is the number of parts of $\lambda$ equal to $1$.
\end{proposition}

\begin{proof}
The number of north steps is the number of gaps, and the path has $F+1$ steps in total, namely $\lambda_1$ east steps
and $\ell(\lambda)$ north steps, so $F=\lambda_1+\ell(\lambda)-1=h_{(1,1)}(\lambda)$. The gap $x$ corresponding to row
$i$ is the index of its north step; the steps after it consist of $\lambda_1-\lambda_i$ east steps and $i-1$ north
steps. So $x=F-(\lambda_1-\lambda_i)-(i-1)=\lambda_i+\ell(\lambda)-i=h_{(i,1)}(\lambda)$. Finally, if $S$ is a
semigroup with multiplicity $m$, the path begins with one east step (the element $0$) followed by $m-1$ north steps
(the gaps $1,\dots,m-1$). These give exactly $m-1$ rows of length $1$, and all later rows are longer.
\end{proof}

\begin{example}\label{ex:KN}
The numerical set $S=\{0,3,4,7,9,10,12,\to\}$ has gaps $\{1,2,5,6,8,11\}$, so $g=6$, $F=11$ and least positive
element $3$. It is not a semigroup, since $3+3=6$ is a gap. Its
Keith--Nath path is $E\,N\,N\,E\,E\,N\,N\,E\,N\,E\,E\,N$, and $\lambda(S)=(6,4,3,3,1,1)$. Indeed $\lambda$ has $6$
parts, $h_{(1,1)}=6+6-1=11$, exactly two parts equal to $1$ (one less than the least positive element), and
first-column hook lengths $11,8,6,5,2,1$, which are the gaps. (This is the example of \cite[Figure 1]{BNS}.)
\end{example}

Following \cite{BNS}, a partition $\lambda$ corresponds to a numerical
semigroup if $S^{KN}(\lambda)$ is a numerical semigroup. We write $\PG(n,g)$, $\PM(n,m)$ and $\PF(n,f)$ for the
numbers of such $\lambda\vdash n$ with, respectively, exactly $g$ parts, exactly $m-1$ parts equal to $1$, and largest
hook length $f$. By Proposition~\ref{prop:KNdict}, these count partitions of $n$ corresponding to numerical semigroups
of genus $g$, multiplicity $m$ and Frobenius number $f$. Finally, let $S'(n)$ be the number of partitions of size at
most $n$ corresponding to numerical semigroups, and $P(n)=\sum_{j\le n}p(j)$; Constantin, Houston-Edwards and Kaplan
conjectured that $S'(n)/P(n)\to0$ \cite{CHK}.

\section{Threshold graphs, creation sequences and numerical sets}\label{sec:cs}

\subsection{Threshold graphs and their creation sequences}

\begin{definition}\label{def:threshold}
A graph $G=(V,E)$ is a \emph{threshold graph} if there are weights $\mathrm{wt}:V\to\mathbb R$ and a threshold
$t\in\mathbb R$ such that, for distinct $u,v\in V$, $\{u,v\}\in E$ if and only if $\mathrm{wt}(u)+\mathrm{wt}(v)\ge t$.
\end{definition}

Equivalently \cite{HIS,MP}, $G$ is a threshold graph if and only if it can be built from a single vertex $v_0$ by
adding vertices $v_1,\dots,v_{N-1}$ one at a time, each new vertex being either \emph{isolated} (adjacent to no
earlier vertex) or \emph{dominating} (adjacent to every earlier vertex).

\begin{definition}\label{def:graphcs}
The \emph{creation sequence} of such a construction is the word $s=s_1s_2\cdots s_{N-1}$ with $s_i=1$ if $v_i$ is
dominating and $s_i=0$ if $v_i$ is isolated. The first vertex $v_0$ involves no choice and receives no symbol.
Conversely, every binary word $s$ of length $N-1$ determines a threshold graph $G(s)$ on $\{v_0,\dots,v_{N-1}\}$, in
which, for $i<j$, $v_i$ is adjacent to $v_j$ if and only if  $s_j=1$. 
\end{definition}

\begin{lemma}[\cite{MP}]\label{lem:graphword}
The map $s\mapsto G(s)$ induces a bijection between binary words of length $N-1$ and isomorphism classes of threshold
graphs on $N$ vertices. In particular there are $2^{N-1}$ threshold graphs on $N$ vertices up to isomorphism.
\end{lemma}

We write $s(G)$ for the creation sequence of a threshold graph $G$ given by Lemma~\ref{lem:graphword}.

\begin{lemma}\label{lem:isolated}
Let $N\ge2$. A threshold graph $G$ on $N$ vertices has an isolated vertex if and only if $s(G)$ ends in $0$.
\end{lemma}

\begin{proof}
If $s(G)$ ends in $0$, the last vertex is adjacent to no earlier vertex and no vertex is added after it, so it is
isolated in $G$. If $s(G)$ ends in $1$, the last vertex is adjacent to all other vertices, so no vertex is isolated.
\end{proof}

\subsection{From threshold graphs to partitions and numerical sets}

\begin{definition}\label{def:graphpartition}
Let $G$ be a threshold graph whose creation sequence $w=s(G)$ ends in $0$. Read $1\,w$ (a leading $1$ for $v_0$,
followed by $w$) as a lattice path: an east step for each $1$ and a north step for each $0$. This path is the profile
of a Young diagram, and we call the corresponding partition $\lambda_G$. Its rows correspond to the isolated vertices
$v_i$ ($i\ge1$), and the row of $v_i$ has length equal to the number of vertices $v_j$ with $j<i$ that are $v_0$ or
dominating. The numerical set of $G$ is
\[
S_G=\{0\}\cup\{i\ge1:\ v_i \text{ is dominating}\}\cup\{N,N+1,\dots\}.
\]
Conversely, for a numerical set $S\ne\N_0$ with Frobenius number $F$, let $G(S)$ be the threshold graph on
$v_0,\dots,v_F$ in which $v_i$ ($1\le i\le F$) is dominating if $i\in S$ and isolated if $i\notin S$.
\end{definition}

\begin{figure}[ht]
\centering
\begin{tikzpicture}[scale=1.25,
  v/.style={circle,fill=black,inner sep=1.6pt},
  lab/.style={font=\small}]
\foreach \i in {0,...,7}{
  \coordinate (p\i) at ({90-45*\i}:1.6);
}
\foreach \j/\k in {4/0,4/1,4/2,4/3,5/0,5/1,5/2,5/3,5/4,6/0,6/1,6/2,6/3,6/4,6/5}{
  \draw[gray] (p\j)--(p\k);
}
\foreach \i in {0,...,7}{
  \node[v] at (p\i){};
  \node[lab] at ({90-45*\i}:1.95) {$v_{\i}$};
}
\end{tikzpicture}
\caption{The threshold graph $G(S)$ for $S=\langle4,5,6\rangle=\{0,4,5,6,8,\dots\}$. Its creation sequence is
$s=0001110$: $v_1,v_2,v_3,v_7$ are isolated and $v_4,v_5,v_6$ are dominating. It has $4+5+6=15$ edges, $v_7$ is an
isolated vertex, and $\lambda_{G}=(4,1,1,1)$.}
\label{fig:threshold}
\end{figure}
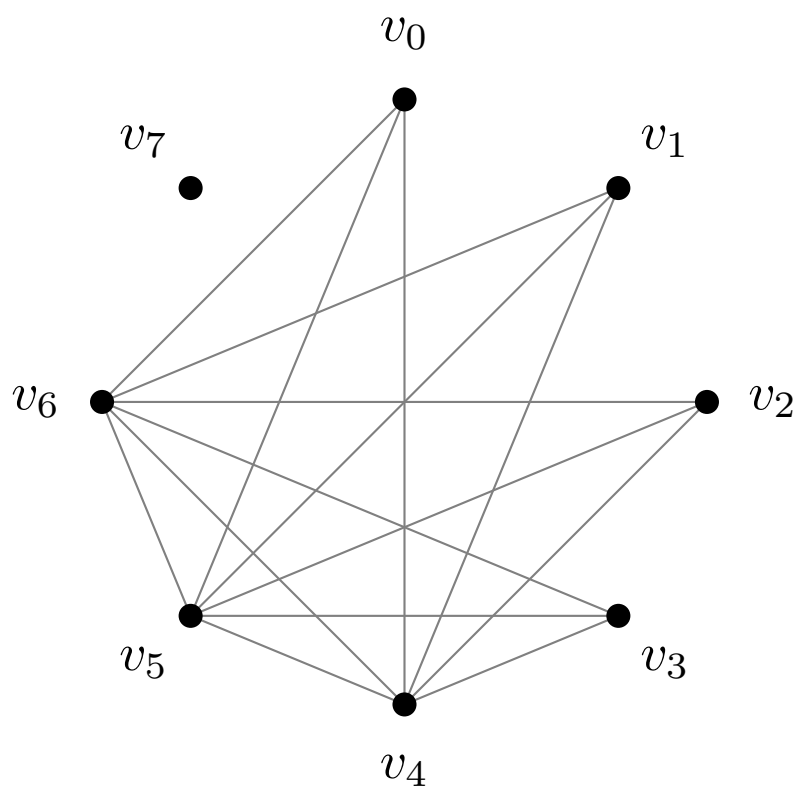

The partition $\lambda_G$ is not the degree partition of $G$, which is the partition usually attached to a threshold
graph \cite{MP}. It is the lattice-path reading of the creation sequence, and it is chosen so that $S_G$ is the
Keith--Nath set of $\lambda_G$ (Theorem~\ref{thm:graphbij}). Figure~\ref{fig:threshold} shows an example.

\subsection{The creation sequence of a partition}

A partition determines the creation sequence of a threshold graph. Most of our proofs are phrased in
terms of partitions whose smallest part is $1$, which are the only ones that can correspond to numerical semigroups, so
we record this direction explicitly.

\begin{definition}\label{def:cs}
For a partition $\lambda=(n_1\ge n_2\ge\dots\ge n_k)$ define its creation sequence
\[
\cs(\lambda)=0\,1^{\,n_{k-1}-n_k}\,0\,1^{\,n_{k-2}-n_{k-1}}\,0\,\cdots\,0\,1^{\,n_1-n_2}\,0 .
\]
This is a binary word with exactly $k$ zeros and $n_1-n_k$ ones, and it always ends in $0$. For a binary word
$w=c_1c_2\cdots c_F$ ending in $0$ put
\[
S_w=\{0\}\cup\{i:\ c_i=1\}\cup\{F+1,F+2,\dots\},
\]
and write $S_\lambda=S_{\cs(\lambda)}$. Thus the letter at position $i$ records whether $i\in S_\lambda$. In the
language of Definition~\ref{def:graphcs}, $\cs(\lambda)$ is the creation sequence of the threshold graph
$G(\cs(\lambda))$: its vertex $v_i$ is dominating exactly when $i\in S_\lambda$.
\end{definition}

The word $\cs(\lambda)$ only sees the differences between consecutive parts and the number of parts. Therefore, 
$\cs(\lambda)=\cs(\lambda+c(1^k))$ for every shift $c$ that keeps the parts positive. For instance $(2,1)$, $(4,3)$
and $(5,4)$ all have creation sequence $010$.

\begin{lemma}\label{lem:bij}
\begin{enumerate}
\item The map $w\mapsto S_w$ is a bijection from binary words ending in $0$ onto numerical sets different from
$\N_0$. The gaps of $S_w$ are the positions of the zeros, so $F(S_w)=|w|$ and $g(S_w)$ is the number of zeros.
\item For a nonempty partition $\lambda$ with smallest part $n_k$, put
\[
w^{KN}(\lambda)=1^{\,n_k-1}\,\cs(\lambda).
\]
Then $S^{KN}(\lambda)=S_{w^{KN}(\lambda)}$. The map $\lambda\mapsto w^{KN}(\lambda)$ is a bijection from the nonempty
partitions onto the binary words ending in $0$. It restricts to a bijection $\lambda\mapsto\cs(\lambda)$ from partitions with smallest
part $1$ onto the words that begin and end with $0$; for such $\lambda$, $S_\lambda=S^{KN}(\lambda)$.
\item For an arbitrary partition, $S_\lambda=S^{KN}(\hat\lambda)$ with $\hat\lambda=(n_1-n_k+1,\dots,n_{k-1}-n_k+1,1)$.
If $n_k\ge2$, then $1\in S^{KN}(\lambda)$, and $S^{KN}(\lambda)$ is not a numerical semigroup.
\end{enumerate}
\end{lemma}

\begin{proof}
(1) A numerical set $S\ne\N_0$ is determined by its gap set, which is finite and nonempty with maximum $F(S)$.
Recording membership of $1,\dots,F(S)$ gives a word of length $F(S)$ whose last letter is $0$. This construction is
inverse to $w\mapsto S_w$. The statements about gaps, $F$ and $g$ are immediate from the definition of $S_w$.

(2) Follow the Keith--Nath path of $\lambda$ from the bottom-left corner. The bottom row has length $n_k$, so the
path begins with $n_k$ east steps; the first of them is the element $0$. The path then takes one north step for
each row. Between the north step of a row of length $n_{i+1}$ and that of the next row up, of length $n_i$, it takes
exactly $n_i-n_{i+1}$ east steps. Reading the steps after the first one, with $1$ for east and $0$ for north, gives
$1^{n_k-1}$ followed by the letters of $\cs(\lambda)$, that is, the word $w^{KN}(\lambda)$. The last step is the north
step of the top row, which is the gap $F$. Thus, $S^{KN}(\lambda)=S_{w^{KN}(\lambda)}$. Conversely, a word ending in $0$
with $k$ zeros can be written uniquely as $1^{e}\,0\,1^{d_{k-1}}\,0\cdots1^{d_1}\,0$ with $e,d_i\ge0$, and it equals
$w^{KN}(\lambda)$ for exactly one partition, namely the one with $k$ parts, $n_k=e+1$ and $n_i-n_{i+1}=d_i$. The word
begins with $0$ exactly when $e=0$, that is, $n_k=1$, and then $w^{KN}(\lambda)=\cs(\lambda)$.

(3) The partitions $\lambda$ and $\hat\lambda$ have the same number of parts and the same differences between
consecutive parts, so $\cs(\lambda)=\cs(\hat\lambda)$, and $\hat\lambda$ has smallest part $1$. If $n_k\ge2$, the
Keith--Nath path of $\lambda$ begins with $n_k\ge2$ east steps, so $0,1\in S^{KN}(\lambda)$. A numerical semigroup
containing $1$ is $\N_0$, but $S^{KN}(\lambda)\neq\N_0$, since $\lambda\ne\emptyset$.
\end{proof}

\begin{proposition}\label{prop:dict}
Let $\lambda=(n_1\ge\dots\ge n_k)$ and let $a=\#\{i:\ n_i=n_k\}$ be the number of copies of its smallest part. Then
\[
F(S_\lambda)=|\cs(\lambda)|=k+n_1-n_k=h_{(1,1)}(\hat\lambda),\qquad g(S_\lambda)=k=\ell(\lambda),\qquad
m(S_\lambda)=a+1 .
\]
If $\lambda$ has at least two distinct parts, then $\cs(\lambda)$ begins with exactly $a$ zeros followed by a $1$.
Every maximal block of zeros other than the first, say of length $L$, corresponds to a part value larger than $n_k$
that occurs exactly $L$ times, and every maximal block of ones corresponds to a difference $n_i-n_{i+1}>0$.
\end{proposition}

\begin{proof}
The length and the number of zeros were computed in Definition~\ref{def:cs}, and
$h_{(1,1)}(\hat\lambda)=\hat\lambda_1+k-1=n_1-n_k+k$. In $\cs(\lambda)$ the letters between the $i$-th and
$(i+1)$-st zeros (counting from the left) are $1^{n_{k-i}-n_{k-i+1}}$. The first $a-1$ of these differences vanish
and the $a$-th does not. Therefore the word starts with $0^a1$, and its first $1$ is at position $a+1=m(S_\lambda)$.
More generally, a maximal run of $L$ zeros after the first run means that $L-1$ consecutive differences vanish, with
nonzero differences on both sides (or the end of the word). This happens exactly when some part value occurs $L$
times, and that value exceeds $n_k$ because a $1$ precedes the run. A maximal run of ones is a single nonzero
difference.
\end{proof}

\subsection{The correspondence and its dictionary}

\begin{theorem}\label{thm:graphbij}
Let $N\ge2$. The maps
\[
G\longmapsto s(G)\longmapsto\lambda_G\longmapsto S_G
\]
restrict to bijections between the following four sets:
\begin{enumerate}
\item[(a)] threshold graphs on $N$ vertices with at least one isolated vertex, up to isomorphism;
\item[(b)] binary words of length $N-1$ ending in $0$;
\item[(c)] nonempty partitions with largest hook length $h_{(1,1)}=N-1$;
\item[(d)] numerical sets with Frobenius number $N-1$.
\end{enumerate}
Moreover $s(G)=w^{KN}(\lambda_G)$, $S_G=S_{s(G)}=S^{KN}(\lambda_G)$, and the inverse of $G\mapsto S_G$ is
$S\mapsto G(S)$. The smallest part of $\lambda_G$ is $1$ if and only if $v_1$ is isolated,    that is, $1\notin S_G$, and in
that case $s(G)=\cs(\lambda_G)$.
\end{theorem}

\begin{proof}
 By Lemma~\ref{lem:graphword}, $G\mapsto s(G)$ is a bijection from threshold graphs on $N$
vertices onto words of length $N-1$, and by Lemma~\ref{lem:isolated} the graphs with an isolated vertex correspond
to the words ending in $0$, which establishes the bijection between (a) and (b).

By Lemma~\ref{lem:bij}-(2), the Keith--Nath path of a nonempty partition
$\lambda$ is $1\,w^{KN}(\lambda)$. By Definition~\ref{def:graphpartition}, the path of $\lambda_G$ is $1\,s(G)$. Therefore, 
$w^{KN}(\lambda_G)=s(G)$, and Lemma~\ref{lem:bij}-(2) shows that $w\mapsto\lambda$ with $w^{KN}(\lambda)=w$ is a
bijection from words ending in $0$ onto nonempty partitions. The path has $\lambda_1$ east steps and $\ell(\lambda)$
north steps, so $N=\lambda_1+\ell(\lambda)=h_{(1,1)}(\lambda)+1$. The statement about the smallest part is the last
assertion of Lemma~\ref{lem:bij}-(2), which gives the bijection between (b) and (c).

By Lemma~\ref{lem:bij}-(1), $w\mapsto S_w$ is a bijection onto numerical sets with
$F(S_w)=|w|=N-1$. By Definition~\ref{def:graphpartition}, $S_G=S_{s(G)}$, and $S_{s(G)}=S^{KN}(\lambda_G)$ by
Lemma~\ref{lem:bij}-(2). Finally, $G(S)$ has creation sequence equal to the word of $S$, so $S\mapsto G(S)$ inverts
$G\mapsto S_G$, there we have the bijection between (c) and (d).
\end{proof}

\begin{proposition}\label{prop:notconverse}
\begin{enumerate}
\item Every numerical semigroup $S\ne\N_0$ has a creation sequence: there is a unique threshold graph $G=G(S)$ with an
isolated vertex such that $S_G=S$. Its creation sequence is $s(G)=c_1\cdots c_{F(S)}$, with $c_i=1$ if and only if
$i\in S$.
\item The converse fails. For every $N\ge3$ there are threshold graphs on $N$ vertices with an isolated vertex whose
numerical set $S_G$ is not a numerical semigroup, and for $N\ge2$ there are threshold graphs whose creation sequence
ends in $1$, which do not determine a numerical set with Frobenius number $N-1$.
\end{enumerate}
\end{proposition}

\begin{proof}
(1) This is the bijection between (d) and (a) of Theorem~\ref{thm:graphbij}, restricted to numerical semigroups.

(2) For $N\ge3$, let $s=0\,1\,0^{N-3}$, that is, $v_1$ is isolated, $v_2$ is dominating and $v_3,\dots,v_{N-1}$ are
isolated. Then $S_G=\{0,2\}\cup[N,\infty)$. If $N\ge5$, then $2+2=4\le N-1$ is a gap, so $S_G$ is not a semigroup.
For $N=3$ and $N=4$ use instead $s=1\,0$ and $s=1\,0\,0$, then $1\in S_G$, but $S_G\ne\N_0$, so $S_G$ is not a
semigroup. Finally, the graph with creation sequence $1^{N-1}$ (the complete graph $K_N$) has no isolated vertex, by
Lemma~\ref{lem:isolated}.
\end{proof}

Thus creation sequences (equivalently, threshold graphs with an isolated vertex) are in bijection with
\emph{numerical sets}, and the numerical semigroups form a proper subclass. By Corollary~\ref{cor:graphcount}, the
proportion of threshold graphs that yield semigroups is the proportion of numerical sets with Frobenius number $N-1$ that
are semigroups. There are $2^{N-2}$ numerical sets with Frobenius number $N-1$, but only of the order of
$2^{N/2}$ numerical semigroups with that Frobenius number \cite{Backelin}, so this proportion tends to $0$ as
$N\to\infty$.

\begin{example}\label{ex:fivevertices}
Table~\ref{tab:five} lists all $2^{3}=8$ threshold graphs on $N=5$ vertices that have an isolated vertex, that is, all
creation sequences $s=c_1c_2c_3c_4$ with $c_4=0$. In each case $S_G=\{0\}\cup\{i:c_i=1\}\cup[5,\infty)$. Only two of
the eight creation sequences give numerical semigroups: $0000$, giving $\{0,5,6,\dots\}=\langle5,6,7,8,9\rangle$, and
$0010$, giving $\{0,3,5,6,\dots\}=\langle3,5,7\rangle$. For example, the creation sequence $0100$ builds the
threshold graph in which $v_2$ is adjacent to $v_0,v_1$ and $v_3,v_4$ are isolated. Its partition is $(2,2,1)$ and its
numerical set $\{0,2,5,6,\dots\}$ is not closed under addition, because $v_2$ is dominating while $v_{2+2}=v_4$ is
isolated.
\begin{table}[ht]
\centering
\begin{tabular}{cccl}
\toprule
$s(G)$ & $\lambda_G$ & $S_G\cap[0,4]$ & numerical semigroup?\\
\midrule
$0000$ & $(1,1,1,1)$ & $\{0\}$ & yes\\
$0010$ & $(2,1,1)$ & $\{0,3\}$ & yes\\
$0100$ & $(2,2,1)$ & $\{0,2\}$ & no: $2+2=4$\\
$0110$ & $(3,1)$ & $\{0,2,3\}$ & no: $2+2=4$\\
$1000$ & $(2,2,2)$ & $\{0,1\}$ & no: $1+1=2$\\
$1010$ & $(3,2)$ & $\{0,1,3\}$ & no: $1+1=2$\\
$1100$ & $(3,3)$ & $\{0,1,2\}$ & no: $1+2=3$\\
$1110$ & $(4)$ & $\{0,1,2,3\}$ & no: $1+3=4$\\
\bottomrule
\end{tabular}
\caption{The eight threshold graphs on five vertices with an isolated vertex. Only two of their creation sequences
give numerical semigroups.}
\label{tab:five}
\end{table}
\end{example}

\begin{proposition}\label{prop:graphdict}
Let $G$ be a threshold graph on $N\ge2$ vertices with an isolated vertex, and put $S=S_G$, $\lambda=\lambda_G$ and
$g=g(S)$. Then the following holds.
\begin{enumerate}
\item $F(S)=N-1=h_{(1,1)}(\lambda)$;
\item $g(S)$ is the number of isolated vertices $v_i$ with $i\ge1$ (isolated at the moment they are added), and equals
$\ell(\lambda)$;
\item $m(S)$ is the least $i\ge1$ such that $v_i$ is dominating ($m(S)=N$ if there is none), and equals $1$ plus the
number of parts of $\lambda$ equal to $1$;
\item the number of edges of $G$ is $\displaystyle\sum_{x\in S,\ 0<x<F(S)}x$;
\item $|\lambda|$ is the number of non-adjacent pairs $\{v_i,v_j\}$ with $i<j$ such that $v_i$ is $v_0$ or dominating
(and $v_j$ is isolated). Equivalently, $|\lambda|=\#\{\text{non-edges of }G\}-\binom g2$;
\item $S$ is a numerical semigroup if and only if, for all $i,j\ge1$ with $i+j\le N-1$, the vertex $v_{i+j}$ is
dominating whenever $v_i$ and $v_j$ are;
\item the complement $\overline G$ is a threshold graph whose creation sequence is obtained from $s(G)$ by
interchanging $0$ and $1$.
\end{enumerate}
\end{proposition}

\begin{proof}
(1)--(3) follow from Theorem~\ref{thm:graphbij}, Lemma~\ref{lem:bij}(1) and Proposition~\ref{prop:dict}, since
the gaps of $S$ in $[1,N-1]$ are the indices of the isolated vertices and the elements are the indices of the
dominating ones.

(4) Every edge $\{v_i,v_j\}$ with $i<j$ is created when $v_j$ is added, and $v_j$ is then dominating, so $j\in S$. A
dominating vertex $v_x$ creates exactly $x$ edges, one to each of $v_0,\dots,v_{x-1}$.

(5) A pair $\{v_i,v_j\}$ with $i<j$ is a non-edge exactly when $v_j$ is isolated, that is,  $j$ is a gap. There are
$\binom g2$ such pairs in which $v_i$ is also isolated ($i\ge1$). The remaining non-edges are counted, for each gap
$j$, by the number of elements of $S$ smaller than $j$, which is the length of the row of $\lambda$ belonging to $j$
(Definition~\ref{def:graphpartition}). Summing over the gaps gives $|\lambda|$.

(6) This is closure under addition for sums of nonzero elements up to $F(S)=N-1$. Sums involving $0$, or sums
exceeding $F(S)$, lie in $S$ automatically.

(7) For $i<j$, $v_i$ and $v_j$ are adjacent in $\overline G$ exactly when $s_j=0$. So $\overline G=G(\bar s)$, where
$\bar s$ is the complementary word.
\end{proof}

\begin{corollary}\label{cor:graphcount}
The number of threshold graphs on $N$ vertices, up to isomorphism, that have an isolated vertex and whose numerical
set is a numerical semigroup equals the number of numerical semigroups with Frobenius number $N-1$.
\end{corollary}

\begin{proof}
This follows from Theorem~\ref{thm:graphbij}.
\end{proof}

\begin{remark}\label{rem:complement}
By Proposition~\ref{prop:graphdict}-(7), complementation of threshold graphs interchanges $0$ and $1$ in the creation
sequence. It leaves the class (a) whenever $s(G)$ ends in $0$, since then $s(\overline G)$ ends in $1$. Conjugation
of partitions corresponds instead to complementing and reversing the word (Lemma~\ref{lem:conj}).
\end{remark}

\begin{definition}\label{def:KNset}
For a partition $\kappa=(\kappa_1\ge\dots\ge\kappa_r)$ let
\[
w(\kappa)=1^{\kappa_r}\,0\,1^{\kappa_{r-1}-\kappa_r}\,0\cdots1^{\kappa_1-\kappa_2}\,0,
\]
a word of length $\kappa_1+r$ whose positions we index by $0,1,\dots,\kappa_1+r-1$. Let $\Sigma(\kappa)\subseteq\N_0$ be
the set of positions of its ones together with all integers $\ge\kappa_1+r$. We put $\Sigma(\emptyset)=\N_0$, and
we treat negative integers as never lying in $\Sigma(\kappa)$.
\end{definition}

\begin{lemma}\label{lem:Sigma}
Let $\kappa$ be a nonempty partition with $r$ parts. Then the following holds.
\begin{enumerate}
\item $\Sigma(\kappa)=S^{KN}(\kappa)$. In particular $0\in\Sigma(\kappa)$ and $[\kappa_1+r,\infty)\subseteq\Sigma(\kappa)$.
\item The zeros of $w(\kappa)$ are exactly at the positions $z_t=\kappa_t+r-t$ for $t=1,\dots,r$, and
$z_1>z_2>\dots>z_r$.
\item $\Sigma(\kappa)+\Sigma(\kappa)\supseteq\Sigma(\kappa)$; in particular every integer $\ge\kappa_1+r$ lies in
$\Sigma(\kappa)+\Sigma(\kappa)$.
\end{enumerate}
\end{lemma}

\begin{proof}
(1) Reading $w(\kappa)$ from position $0$ is exactly the Keith--Nath path of $\kappa$ (east for $1$, north for $0$),
which starts with $\kappa_r\ge1$ east steps, so position $0$ carries a $1$. 

(2) The zero that ends the block
$1^{\kappa_t-\kappa_{t+1}}$ is preceded by $\kappa_t$ ones in total (the lengths of the blocks telescope) and by
$r-t$ zeros. Its position is therefore $\kappa_t+r-t$, and these positions strictly decrease in $t$ because
$\kappa_t\ge\kappa_{t+1}$. 

(3) This follows from $0\in\Sigma(\kappa)$ and (1).
\end{proof}

%======================================================================
\section{Structural criteria}\label{sec:struct}
%======================================================================

Throughout this section $w$ is a binary word ending in $0$, and $S=S_w$.

\begin{theorem}\label{thm:short}
If the first $1$ of $w$ is at position $m$ and $|w|<2m$, then $S_w$ is a numerical semigroup. In partition language:
if $\lambda$ has smallest part repeated $a$ times and $k+n_1-n_k<2a+2$, then $S_\lambda$ is a numerical semigroup.
\end{theorem}

\begin{proof}
Every nonzero element of $S$ is at least $m$, so the sum of two nonzero elements is at least $2m>|w|=F(S)$ and
therefore lies in $S$. Sums involving $0$ are trivially in $S$. The second statement is the first combined with
Proposition~\ref{prop:dict}.
\end{proof}

\begin{theorem}\label{thm:blocks}
Let $S_w$ be a numerical semigroup with multiplicity $m$, so that $w=0^{m-1}1\cdots$. Then the following holds.
\begin{enumerate}
\item every maximal block of ones in $w$ has length at most $m-1$; and
\item every maximal block of zeros in $w$ other than the first has length at most $m-1$.
\end{enumerate}
\end{theorem}

\begin{proof}
(1) Let the ones occupy the positions $u,u+1,\dots,u+L-1$ with $L\ge m$. The position $u+L$ carries a $0$: the
block is maximal, and it cannot end the word because $w$ ends in $0$. Every integer $y\ge u$ can be written as
$y=u+i+qm$ with $0\le i\le m-1$ and $q\ge0$. Since $u+i\in S$ and $m\in S$, closure gives $y\in S$. In particular
$u+L\in S$, a contradiction.

(2) Let the zeros occupy the positions $u,\dots,u+L-1$ with $L\ge m$, and suppose this block is not the first one.
Then position $u-1$ carries a $1$, and $u-1\ge m$ because the first $1$ is at $m$. Now
$(u-1)+m\in[u,u+m-1]\subseteq[u,u+L-1]$ is a gap, but it is the sum of two elements of $S$, a contradiction.
\end{proof}

\begin{corollary}\label{cor:blocks}
Let $\lambda=(n_1\ge\dots\ge n_k)$, and let $a$ be the number of copies of its smallest part. If $S_\lambda$ is a
numerical semigroup, then the following holds.
\begin{enumerate}
\item consecutive distinct parts differ by at most $a$: $n_i-n_{i+1}\le a$ for all $i$; and
\item no part value occurs more often than the smallest part: every part value occurs at most $a$ times.
\end{enumerate}
In particular (\emph{long first block}): if the next part value above $n_k$ exceeds $n_k$ by at least $a+1$, then
$S_\lambda$ is not a numerical semigroup.
\end{corollary}

\begin{proof}
By Proposition~\ref{prop:dict}, $m(S_\lambda)=a+1$, the blocks of ones are the nonzero differences, and the
non-initial blocks of zeros are the multiplicities of the part values above $n_k$. Apply
Theorem~\ref{thm:blocks}.
\end{proof}

\begin{remark}
Corollary~\ref{cor:blocks}-(1) is the partition form of the observation in \cite[Section 2.3]{BNS} that ``the gap
between two parts is no larger than $m_\lambda$''.
\end{remark}

\begin{theorem}\label{thm:unique}
Suppose the smallest part of $\lambda$ occurs exactly once, that is, $a=1$, and $k\ge2$. Then $S_\lambda$ is a
numerical semigroup if and only if $\cs(\lambda)=(01)^{k-1}0$, that is, if and only if the parts of $\lambda$ are
distinct consecutive integers $\lambda=(c+k-1,c+k-2,\dots,c)$. In this case $S_\lambda=\langle2,2k+1\rangle$.
\end{theorem}

\begin{proof}
If $\cs(\lambda)=(01)^{k-1}0$, the ones sit at the even positions $2,4,\dots,2k-2$ and the zeros at the odd positions
$1,3,\dots,2k-1$. Thus $S_\lambda=\{0,2,4,\dots,2k-2\}\cup[2k,\infty)$, which consists of all even numbers and all
odd numbers $\ge2k+1$; this is $\langle2,2k+1\rangle$.

Conversely, let $S_\lambda$ be a semigroup. Its multiplicity is $a+1=2$, so by Theorem~\ref{thm:blocks} every block
of ones and every non-initial block of zeros has length $1$. The first block of zeros has length $a=1$. Thus $w$
alternates, $w=(01)^{j}0$, and the number of zeros gives $j=k-1$. Finally, $w=(01)^{k-1}0$ means that all
differences equal $1$ and no part repeats.
\end{proof}

\begin{theorem}\label{thm:equal}
If $\lambda=(c^k)$, then $\cs(\lambda)=0^k$ and $S_\lambda=\{0\}\cup[k+1,\infty)$ is a numerical semigroup.
\end{theorem}

\begin{proof}
All differences vanish, so $\cs(\lambda)=0^k$. The nonzero elements of $S_\lambda$ are at least $k+1$, so the sum of
two of them is at least $2k+2>k=F(S_\lambda)$.
\end{proof}

\begin{theorem}\label{thm:twosizes}
Let $\lambda=(A^p,B^q)$ with $A>B\ge1$, $p,q\ge1$, and put $d=A-B$. Then $S_\lambda$ is a numerical semigroup if and
only if
\[
d+p\le q+1 .
\]
\end{theorem}

\begin{proof}
Here $\cs(\lambda)=0^q1^d0^p$, so $m=q+1$, the elements in $(0,F]$ are $q+1,\dots,q+d$, the gaps are
$1,\dots,q$ and $q+d+1,\dots,q+d+p$, and $F=q+d+p$. If $d+p\le q+1$, then $F\le2q+1<2m$ and
Theorem~\ref{thm:short} applies. Suppose instead that $d+p\ge q+2$. If $d\ge q+1$, then $S_\lambda$ is not a semigroup
by Theorem~\ref{thm:blocks}-(1). If $d\le q$, then $2m=2q+2$ satisfies $q+d+1\le2q+2$ (because $d\le q+1$) and
$2q+2\le q+d+p=F$ (because $d+p\ge q+2$). So $2m$ is a gap while $m\in S_\lambda$, and $S_\lambda$ is not a
semigroup.
\end{proof}

\begin{corollary}\label{cor:hook}
Let $\lambda=(n_1,n_2^{\,k-1})$ with $n_1>n_2$ and $d=n_1-n_2$. Then $S_\lambda$ is a numerical semigroup if and only
if $d\le k-1$. In particular $(n,1^{j})$ corresponds to a numerical semigroup if and only if $j\ge n-1$, and
$\PF(n,n)=\lceil n/2\rceil$ \cite[Theorem 6]{BNS}.
\end{corollary}

\begin{proof}
Take $p=1$ and $q=k-1$ in Theorem~\ref{thm:twosizes}. The partitions of $n$ with largest hook $n$ are exactly the
hooks $(n-j,1^{j})$ with $0\le j\le n-1$. For $1\le j\le n-2$, Theorem~\ref{thm:twosizes} with $d=n-j-1$, $p=1$,
$q=j$ shows that $(n-j,1^j)$ corresponds to a semigroup if and only if $n-j\le j+1$, that is, $j\ge\lfloor n/2\rfloor$.
For $j=n-1$, $(1^n)$ corresponds to a semigroup (Theorem~\ref{thm:equal}). For $j=0$ and $n\ge2$, $(n)$ has no part
$1$. In all cases the condition is $j\ge\lfloor n/2\rfloor$, which holds for $n-\lfloor n/2\rfloor=\lceil n/2\rceil$
values of $j$.
\end{proof}

The following criterion decides every creation sequence. It is Kunz's description of numerical semigroups of fixed
multiplicity \cite{Kunz,RGS}, restated in the language of words.

\begin{theorem}\label{thm:column}
Let $w=c_1\cdots c_F$ end in $0$ and contain at least one $1$, let $m$ be the position of its first $1$, and set
$c_i=1$ for $i>F$. Say that
position $i$ lies in column $i\bmod m$. Then $S_w$ is a numerical semigroup if and only if
\begin{enumerate}
\item[(i)] $c_i=1$ implies $c_{i+m}=1$ for every $i\ge1$; equivalently, column $0$ contains no zero and each column
$r\in\{1,\dots,m-1\}$ reads $0\cdots01\cdots1$ from top to bottom; and
\item[(ii)] the numbers $k_r$ of zeros in column $r$ satisfy
\begin{align*}
k_i+k_j&\ge k_{i+j} && (1\le i\le j,\ i+j<m),\\
k_i+k_j+1&\ge k_{i+j-m} && (1\le i\le j\le m-1,\ i+j>m).
\end{align*}
\end{enumerate}
In that case $g(S_w)=\sum_rk_r$ and $F(S_w)=\max_r(k_rm+r)-m$.
\end{theorem}

\begin{proof}
Put $w_0=0$ and, for $1\le r\le m-1$, $w_r=k_rm+r$.

Assume that $S=S_w$ is a semigroup. Since $m\in S$, $x\in S$ implies $x+m\in S$, which is (i). Consequently the
column $r$ consists of $k_r$ zeros followed by ones, so the least element of $S$ congruent to $r$ modulo $m$ is
$w_r$, and $x\in S$ if and only if $x\ge w_{x\bmod m}$. For $i,j\ge1$, the element $w_i+w_j\in S$ is congruent to
$i+j$, so $w_i+w_j\ge w_{(i+j)\bmod m}$. For $i+j<m$ this reads $(k_i+k_j)m+i+j\ge k_{i+j}m+i+j$, and for $i+j>m$ it
reads $(k_i+k_j+1)m+(i+j-m)\ge k_{i+j-m}m+(i+j-m)$. This is (ii).

Conversely, assume (i) and (ii). By (i), $x\in S$ if and only if $x\ge w_{x\bmod m}$. Let $x,y\in S$ with residues
$i,j$. If $i=0$, then $x$ is a multiple of $m$ and $x+y\ge y\ge w_j$, so $x+y\in S$; similarly if $j=0$ or $i+j=m$
(then $w_0=0$). Otherwise (ii) gives $x+y\ge w_i+w_j\ge w_{(i+j)\bmod m}$, so $x+y\in S$.

Finally, the gaps are the $k_r$ zeros of each column, and the largest gap is $\max_r(w_r-m)$.
\end{proof}

\begin{corollary}\label{cor:m23}
Let $g\ge1$. Then the following holds.
\begin{enumerate}
\item Exactly one partition with $g$ parts and exactly one part equal to $1$ corresponds to a numerical semigroup,
namely the staircase $(g,g-1,\dots,1)$, with $S=\langle2,2g+1\rangle$.
\item For $g\ge2$, exactly $\lfloor g/3\rfloor+1$ partitions with $g$ parts and exactly two parts equal to $1$
correspond to numerical semigroups. They are the partitions with $\cs(\lambda)$ having $k_1$ zeros in column $1$
and $k_2=g-k_1$ zeros in column $2$ (modulo $3$), where $g/3\le k_1\le(2g+1)/3$; the corresponding semigroup is
$\langle3,3k_1+1,3k_2+2\rangle$.
\end{enumerate}
\end{corollary}

\begin{proof}
(1) This is Theorem~\ref{thm:unique} with smallest part $1$.

(2) Here $m=3$ and the positions $1,2$ are gaps, so $k_1,k_2\ge1$. The only instances of (ii) are
$(i,j)=(1,1)$, giving $2k_1\ge k_2$, and $(i,j)=(2,2)$, giving $2k_2+1\ge k_1$. With $k_2=g-k_1$ these become
$k_1\ge g/3$ and $k_1\le(2g+1)/3$, and they also force $k_1,k_2\ge1$ when $g\ge2$. The number of integers $k_1$
in $[g/3,(2g+1)/3]$ is $\lfloor(2g+1)/3\rfloor-\lceil g/3\rceil+1$. Checking the three residues $g=3q,3q+1,3q+2$, this
equals $q+1=\lfloor g/3\rfloor+1$ in each case. By Theorem~\ref{thm:column}, each such pair $(k_1,k_2)$ gives
exactly one semigroup, whose Apéry elements are $3,3k_1+1,3k_2+2$.
\end{proof}

\begin{remark}
Theorem~\ref{thm:column} shows that no finite list of equality patterns can classify all partitions: every infinite
pattern family contains both semigroups and non-semigroups, and a pattern only determines which columns are
compared. The number of partitions with $g$ parts corresponding to numerical semigroups is the number $n_g$ of
numerical semigroups of genus $g$, which grows like $C\varphi^g$ with $\varphi$ the golden ratio \cite{Zhai}.
\end{remark}

%======================================================================
\section{Classification by equality patterns}\label{sec:patterns}
%======================================================================

In this section we apply results of Section~\ref{sec:struct} to partitions
described by their pattern of equalities and strict inequalities among the parts.
We first record the consequences for partitions with a unique smallest part, then
treat the pattern $n_1>n_2>n_3=n_4$, and finally determine, for every partition with
at most four parts, exactly when $S_\lambda$ is a numerical semigroup
(Theorem~\ref{thm:small} and Table~\ref{tab:small}).

\begin{corollary}\label{cor:unique}
Let $k\ge2$. Then the following holds.
\begin{enumerate}
\item Every partition into distinct consecutive integers gives a numerical semigroup; for example $(n,n-1)$, the
staircase $(k,k-1,\dots,1)$ and $(4,3,2,1)$.
\item If $n_1>n_2>\dots>n_k$ and some $n_i-n_{i+1}\ge2$, then $S_\lambda$ is not a semigroup. In particular $(n,j)$
with $1\le j\le n-2$ never gives a semigroup.
\item If the smallest part is unique and some part is repeated, then $S_\lambda$ is not a semigroup. This covers the
patterns $n_1=n_2>\dots>n_k$ (with $k\ge3$); \ $n_1>\dots>n_r=n_{r+1}>\dots>n_k$ with $r+1<k$; and
$n_1=\dots=n_{k-1}>n_k$ with $k\ge3$.
\end{enumerate}
\end{corollary}

\begin{proof}
In each case the smallest part is unique, so Theorem~\ref{thm:unique} applies. In (1) the parts are distinct
consecutive integers. In (2) and (3) they are not: in (2) some difference is at least $2$, and in (3) some part is
repeated.
\end{proof}

\begin{proposition}\label{prop:3211}
Let $n_1>n_2>n_3=n_4$, $x=n_2-n_3$ and $y=n_1-n_2$. Then $S_\lambda$ is a semigroup if and only if $(x,y)=(1,2)$, in
which case $S_\lambda=\langle3,5\rangle$.
\end{proposition}

\begin{proof}
Here $\cs(\lambda)=0^21^x01^y0$ and $m=3$. So $3,\dots,x+2\in S_\lambda$, $x+3$ is a gap,
$x+4,\dots,x+y+3\in S_\lambda$, and $F=x+y+4$. If $x\ge3$ or $y\ge3$, Theorem~\ref{thm:blocks}-(1) applies. The
remaining cases are $x,y\in\{1,2\}$ which are listed a follows:

(i) $(x,y)=(1,1)$: $S_\lambda=\{0,3,5,7,\dots\}$, and $3+3=6$ is a gap.\\

(ii) $(x,y)=(2,1)$: $S_\lambda=\{0,3,4,6,8,\dots\}$, and $3+4=7$ is a gap.\\

(iii) $(x,y)=(2,2)$: $S_\lambda=\{0,3,4,6,7,9,\dots\}$, and $4+4=8$ is a gap.\\

(iv) $(x,y)=(1,2)$: $S_\lambda=\{0,3,5,6,8,9,\dots\}=\langle3,5\rangle$.\\

\end{proof}

\begin{theorem}\label{thm:small}
For partitions with at most four parts, $S_\lambda$ is a numerical semigroup exactly in the cases listed in
Table~\ref{tab:small}, where $d_i=n_i-n_{i+1}$.
\end{theorem}

\begin{proof}
We go through the rows of Table~\ref{tab:small}.

(i) \emph{All parts equal:} Theorem~\ref{thm:equal}.
\item \emph{Unique smallest part} (the patterns $n_1>n_2$, $n_1=n_2>n_3$, $n_1>n_2>n_3$, $n_1=n_2=n_3>n_4$,
$n_1=n_2>n_3>n_4$, $n_1>n_2=n_3>n_4$, $n_1>n_2>n_3>n_4$): Theorem~\ref{thm:unique}. A semigroup occurs exactly when
the parts are distinct and consecutive, which is possible only in the strictly decreasing patterns.\\

(ii) \emph{Two part sizes} (the patterns $n_1>n_2=n_3$, $n_1>n_2=n_3=n_4$, $n_1=n_2>n_3=n_4$):
Theorem~\ref{thm:twosizes}, with $(p,q)=(1,2),(1,3),(2,2)$ respectively, giving $d_1\le2$, $d_1\le3$ and $d_2\le1$.\\

(iii) \emph{The pattern $n_1>n_2>n_3=n_4$:} Proposition~\ref{prop:3211}. \qedhere

\end{proof}

\begin{table}[ht]
\centering
\begin{tabular}{clll}
\toprule
$k$ & pattern & $S_\lambda$ is a semigroup if and only if & by\\
\midrule
2 & $n_1=n_2$ & always & Thm~\ref{thm:equal}\\
2 & $n_1>n_2$ & $d_1=1$ & Thm~\ref{thm:unique}\\
\midrule
3 & $n_1=n_2=n_3$ & always & Thm~\ref{thm:equal}\\
3 & $n_1=n_2>n_3$ & never & Thm~\ref{thm:unique}\\
3 & $n_1>n_2=n_3$ & $d_1\le2$ & Thm~\ref{thm:twosizes}\\
3 & $n_1>n_2>n_3$ & $d_1=d_2=1$ & Thm~\ref{thm:unique}\\
\midrule
4 & $n_1=n_2=n_3=n_4$ & always & Thm~\ref{thm:equal}\\
4 & $n_1=n_2=n_3>n_4$ & never & Thm~\ref{thm:unique}\\
4 & $n_1=n_2>n_3=n_4$ & $d_2=1$ & Thm~\ref{thm:twosizes}\\
4 & $n_1=n_2>n_3>n_4$ & never & Thm~\ref{thm:unique}\\
4 & $n_1>n_2=n_3=n_4$ & $d_1\le3$ & Thm~\ref{thm:twosizes}\\
4 & $n_1>n_2=n_3>n_4$ & never & Thm~\ref{thm:unique}\\
4 & $n_1>n_2>n_3=n_4$ & $(d_2,d_1)=(1,2)$ & Prop~\ref{prop:3211}\\
4 & $n_1>n_2>n_3>n_4$ & $d_1=d_2=d_3=1$ & Thm~\ref{thm:unique}\\
\bottomrule
\end{tabular}
\caption{Partitions with at most four parts whose creation sequence gives a numerical semigroup.}
\label{tab:small}
\end{table}

\begin{example}\label{ex:patterns}
\begin{enumerate}
\item $(4,3)$: $\cs=010$, $S=\{0,2,4,5,\dots\}=\langle2,5\rangle$.
\item $(5,5,3)$: $\cs=01^20^2$, $S=\{0,2,3,6,\dots\}$; $2+2=4\notin S$.
\item $(6,4,4)$: $\cs=0^21^20$, $S=\{0,3,4,6,\dots\}=\langle3,4\rangle$.
\item $(6,4,3,1)$: $\cs=01^20101^20$, $S=\{0,2,3,5,7,8,10,\dots\}$; $2+2=4\notin S$.
\item $(4,3,2,1)$: $\cs=0101010$, $S=\langle2,9\rangle$.
\item $(4,4,3,2)$: $\cs=01010^2$, $S=\{0,2,4,7,\dots\}$; $2+4=6\notin S$.
\item $(9,7,3,3)$: $\cs=0^21^401^20$, $S=\{0,3,4,5,6,8,9,11,\dots\}$; $3+4=7\notin S$.
\item $(7,5,4,4)$: $\cs=0^2101^20$, $S=\{0,3,5,6,8,\dots\}=\langle3,5\rangle$.
\item $(9,8,8,2)$: $\cs=01^60^210$, $S=\{0,2,\dots,7,10,12,\dots\}$; $2+6=8\notin S$.
\item $(9,9,6,6)$: $\cs=0^21^30^2$, $S=\{0,3,4,5,8,\dots\}$; $3+3=6\notin S$.
\item $(7,7,6,6)$: $\cs=0^210^2$, $S=\{0,3,6,7,\dots\}=\langle3,7,8\rangle$, with gaps $\{1,2,4,5\}$.
\item $(9,9,9,5)$: $\cs=01^40^3$, $S=\{0,2,3,4,5,9,\dots\}$; $3+3=6\notin S$.
\item $(4,1,1,1)$: $\cs=0^31^30$, $S=\{0,4,5,6,8,\dots\}=\langle4,5,6\rangle$.
\item $(7,6,5,4,3,2,1)$: $\cs=(01)^60$, $S=\langle2,15\rangle$.
\item $(11,9,7,5,3,1)$: $\cs=0(1^20)^5$, $S=\{0,2,3,5,6,8,9,11,12,14,15,17,\dots\}$; $2+2=4\notin S$.
\item $(11,11,9,7,5,3,2)$: $\cs=010(1^20)^40$, $S=\{0,2,4,5,7,8,10,11,13,14,17,\dots\}$; $2+4=6\notin S$.
\item $(13,11,7,5,5)$: $\cs=0^21^201^401^20$, $S=\{0,3,4,6,7,8,9,11,12,14,\dots\}$; $4+6=10\notin S$.
\item $(13,11,9,8,8,6,4,2)$: $\cs=01^201^201^20^2101^201^20$; $2+2=4\notin S$.
\item $(5,5,5,5,3)$: $\cs=01^20^4$, $S=\{0,2,3,8,\dots\}$; $2+2=4\notin S$.
\end{enumerate}
\end{example}

\section{Conjugation and symmetry}\label{sec:sym}

In this section we use the full Keith--Nath word:
$
u(\lambda)=c_0c_1\cdots c_F,\qquad c_x=1$ if and only if $x\in S^{KN}(\lambda),
$
where $F=h_{(1,1)}(\lambda)$. It always has $c_0=1$ and $c_F=0$, and $u(\lambda)=1\,\cs(\lambda)$ when the smallest
part of $\lambda$ is $1$. We write $\bar u$ for the complement of $u$ (interchange $0$ and $1$) and $u^{\mathrm{rev}}$
for the reversed word. Recall that a numerical semigroup $S$ is \emph{symmetric} if $x\in S$ if and only if $F(S)-x\notin S$ for
all $x\in\mathbb Z$, and \emph{irreducible} if it is not the intersection of two numerical semigroups properly
containing it.

\begin{lemma}\label{lem:conj}
For every nonempty partition $\lambda$, $u(\lambda')=\bigl(\bar u(\lambda)\bigr)^{\mathrm{rev}}$. Equivalently,
$h_{(1,1)}(\lambda')=h_{(1,1)}(\lambda)=F$ and, for $0\le x\le F$,
$
x\in S^{KN}(\lambda')$ if and only if $F-x\notin S^{KN}(\lambda)$.

\end{lemma}

\begin{proof}
The profile of $\lambda$ is a lattice path of $F+1$ unit steps, going from the lower end of the first column to the
right end of the first row. Transposition in the main diagonal maps the diagram of $\lambda$ onto that of
$\lambda'$, and maps this path onto the profile of $\lambda'$. It exchanges horizontal and vertical steps, and it
traverses the path in the opposite direction: the start of the path of $\lambda$ becomes the end of the path of
$\lambda'$. Therefore the step with index $x$ in $u(\lambda')$ is the complement of the step with index $F-x$ in
$u(\lambda)$. The number of steps is unchanged, so the largest hook lengths agree.
\end{proof}

\begin{theorem}\label{thm:ellbound}
Let $\lambda$ correspond to a numerical semigroup $S$. Then
$
\ell(\lambda)\ \ge\ \lambda_1,
$
and the following are equivalent:

(a) $\ell(\lambda)=\lambda_1$;

(b) $S$ is symmetric;

(c) $\lambda=\lambda'$.

Moreover, if $F(S)$ is even, then $\ell(\lambda)\ge\lambda_1+1$.
\end{theorem}

\begin{proof}
Let $F=F(S)=\lambda_1+\ell(\lambda)-1$ and $g=g(S)=\ell(\lambda)$. Consider the map
$\sigma\colon\{0,\dots,F\}\to\{0,\dots,F\}$, $\sigma(x)=F-x$. Since $\sigma\circ\sigma=\mathrm{id}$, $\sigma$ is an
involution, and the group $\langle\sigma\rangle=\{\mathrm{id},\sigma\}\cong\mathbb{Z}/2\mathbb{Z}$ acts on
$\{0,\dots,F\}$ by $\mathrm{id}\cdot x=x$ and $\sigma\cdot x=F-x$. The orbit of $x$ under this action is
$\{x,\sigma(x)\}=\{x,F-x\}$. The orbits partition $\{0,\dots,F\}$, and they are the pairs $\{x,F-x\}$ of integers
adding up to $F$, together with the fixed point $\{F/2\}$ when $F$ is even (the only $x$ with $\sigma(x)=x$). Each
orbit contains at most one element of $S$: if $x$ and $F-x$ were both in $S$, then $F=x+(F-x)\in S$, a contradiction.

If $F$ is odd, then the action has no fixed points, so there are $(F+1)/2$ orbits, all of size $2$. Thus,
$|S\cap[0,F]|\le(F+1)/2$ and $g=(F+1)-|S\cap[0,F]|\ge(F+1)/2$. If $F$ is even, there are $F/2$ orbits of size $2$ and
the fixed orbit $\{F/2\}$, which contains no element of $S$, since $F/2+F/2=F$. Therefore, $|S\cap[0,F]|\le F/2$ and
$g\ge F/2+1$. Substituting $F=\lambda_1+g-1$ gives $g\ge\lambda_1$ in the first case and $g\ge\lambda_1+1$ in the
second.

Equivalence of (a) and (b): Equality $g=\lambda_1$ forces $F=2g-1$ to be odd and $|S\cap[0,F]|=(F+1)/2$, that is,
every orbit $\{x,F-x\}$ of the action contains exactly one element of $S$. For $0\le x\le F$ this is precisely
symmetry. For $x<0$ or $x>F$ the condition $x\in S$ if and only if $F-x\notin S$ holds trivially: if $x<0$ then
$x\notin S$ and $F-x>F$, so $F-x\in S$, and the case $x>F$ is symmetric. Conversely, symmetry gives exactly one
element of $S$ in each orbit, and the equality follows.

Equivalence of (b) and (c): By Lemma~\ref{lem:conj}, $S$ is symmetric if and only if
$S^{KN}(\lambda')=S^{KN}(\lambda)$, and since the Keith--Nath map is a bijection, this holds if and only if
$\lambda'=\lambda$.
\end{proof}

\begin{remark}
The equivalence of (b) and (c) in Theorem~\ref{thm:ellbound}, that symmetric numerical semigroups correspond
to self-conjugate Young diagrams, is known and plays a central role in \cite{SY21,Yes25}; see also \cite{Tut19} for
Young duals of Arf semigroups. Our contribution here is the inequality $\ell(\lambda)\ge\lambda_1$ and its
consequences below.
\end{remark}

\begin{corollary}\label{cor:irred}
Let $\lambda$ correspond to a numerical semigroup $S$. Then $S$ is irreducible if and only if
$\ell(\lambda)-\lambda_1\in\{0,1\}$. It is symmetric if the difference is $0$, and pseudo-symmetric if it is $1$.
\end{corollary}

\begin{proof}
A numerical semigroup is irreducible if and only if $g(S)=\lceil(F(S)+1)/2\rceil$, and it is then symmetric for $F$
odd and pseudo-symmetric for $F$ even \cite[Chapter 3]{RGS}. With $g=\ell$ and $F=\lambda_1+\ell-1$: for $F$ odd the
condition reads $\ell=\lambda_1$, and for $F$ even it reads $\ell=(\lambda_1+\ell)/2+\tfrac12$, that is, 
$\ell=\lambda_1+1$. A difference of $0$ forces $F$ odd and a difference of $1$ forces $F$ even.
\end{proof}

\begin{corollary}\label{cor:bothconj}
Both $\lambda$ and $\lambda'$ correspond to numerical semigroups if and only if $\lambda=\lambda'$ and the
corresponding semigroup is symmetric.
\end{corollary}

\begin{proof}
Theorem~\ref{thm:ellbound} applied to $\lambda$ and to $\lambda'$ gives $\ell(\lambda)\ge\lambda_1$ and
$\lambda_1=\ell(\lambda')\ge\lambda'_1=\ell(\lambda)$. Thus the equality holds, and $\lambda=\lambda'$ by the theorem.
The converse is clear.
\end{proof}

For the density statement, let $S'(n)$ be the number of partitions of size at most $n$ corresponding to numerical
semigroups, and let $P(n)=\sum_{j\le n}p(j)$, as in \cite{CHK}.

\begin{lemma}\label{lem:equal}
Let $e(n)=\#\{\lambda\vdash n:\ \ell(\lambda)=\lambda_1\}$. Then $e(n)/p(n)\to0$.
\end{lemma}

\begin{proof}
Fix $K\ge1$. For $1\le s\le K$, let $D_s=\{\lambda\vdash n:\ \ell(\lambda)-\lambda_1=2s\}$, 
$E=\{\lambda\vdash n:\ \ell(\lambda)=\lambda_1\}$, and let $B_s=\{\lambda\vdash n:\ \lambda_1-\lambda_2<s\}$ (with
$\lambda_2=0$ if $\ell(\lambda)=1$). For $\lambda\in E\setminus B_s$ with $n\ge2$ define
\[
\phi_s(\lambda)=(\lambda_1-s,\lambda_2,\dots,\lambda_\ell,1^s).
\]
This is a partition because $\lambda_1-s\ge\lambda_2\ge1$. It has $\ell+s$ parts and largest part $\lambda_1-s$, so
$\phi_s(\lambda)\in D_s$. Moreover $\phi_s$ is injective: $\lambda$ is recovered by deleting $s$ parts equal to $1$
and adding $s$ to the first part. The sets $E,D_1,\dots,D_K$ are pairwise disjoint, so
\[
|E|+\sum_{s=1}^K\bigl(|E|-|B_s|\bigr)\le p(n), ~\text{which implies}~ (K+1)\,e(n)\le p(n)+K|B_K| .
\]
By conjugation, $\lambda_1-\lambda_2$ equals the number of parts of $\lambda'$ equal to $1$. Therefore
$$|B_K|=\sum_{t<K}\#\{\mu\vdash n:\ \mu\ \text{has exactly $t$ parts equal to }1\}=\sum_{t<K}\bar p(n-t),$$ and
$\bar p(n-t)=p(n-t)-p(n-t-1)=o(p(n))$ by the Hardy--Ramanujan asymptotic formula, which gives $p(n-1)/p(n)\to 1$
\cite{HR}. Thus $\limsup e(n)/p(n)\le1/(K+1)$ for every $K$.
\end{proof}

\begin{theorem}\label{thm:density}
For every $n$,
\[
\#\{\lambda\vdash n:\ \lambda\ \text{corresponds to a numerical semigroup}\}\ \le\ \frac{p(n)+e(n)}{2},
\]
and consequently $\displaystyle\limsup_{n\to\infty}\frac{S'(n)}{P(n)}\le\frac12$.
\end{theorem}

\begin{proof}
By Theorem~\ref{thm:ellbound}, the partitions in question satisfy $\ell(\lambda)\ge\lambda_1$. Conjugation is an
involution exchanging $\{\ell>\lambda_1\}$ and $\{\ell<\lambda_1\}$, so these two sets have the same size
$(p(n)-e(n))/2$, and
$\#\{\ell\ge\lambda_1\}=(p(n)-e(n))/2+e(n)=(p(n)+e(n))/2$. Summing over $j\le n$ and using $e(j)/p(j)\to0$
(Lemma~\ref{lem:equal}) together with $P(n)\to\infty$ gives the limit statement.
\end{proof}

\begin{remark}
Constantin, Houston-Edwards and Kaplan conjectured that $S'(n)/P(n)\to0$ \cite{CHK}. Theorem~\ref{thm:density} is an
unconditional first step. Numerically the ratio is much smaller: $S'(44)/P(44)\approx0.115$.
\end{remark}

\section{The staircase theorem}\label{sec:staircase}
Let $\delta_g=(g,g-1,\dots,2,1)$ be the staircase partition with $g$ parts. We say $\lambda\subseteq\delta_g$ if
$\lambda_i\le g+1-i$ for every $i$. For a numerical set $S$ with gaps $x_1<x_2<\dots<x_g$, the \emph{weight} of $S$ is
\[
w(S)=\sum_{i=1}^{g}x_i-\frac{g(g+1)}{2}=\sum_{i=1}^{g}(x_i-i)\ \ge0 .
\]
For numerical semigroups arising as Weierstrass semigroups of algebraic curves, this is the classical weight of the
Weierstrass point \cite{EH}.

\begin{lemma}\label{lem:rowgap}
Let $S\ne\N_0$ be a numerical set with gaps $x_1<\dots<x_g$ and Keith--Nath partition $\lambda=\lambda(S)$. The row of
$\lambda$ belonging to the gap $x_i$ is the row $\lambda_{g+1-i}$, and
\[
\lambda_{g+1-i}=x_i-(i-1)\qquad(1\le i\le g).
\]
\end{lemma}

\begin{proof}
Rows correspond to gaps, and they are listed from the top, so the largest gap $x_g$ gives the first row and the
smallest gap $x_1$ gives the last row $\lambda_g$. The length of the row of $x_i$ is the number of elements of $S$ in
$[0,x_i)$. That interval contains $x_i$ integers, of which exactly $i-1$ are gaps.
\end{proof}

\begin{theorem}\label{thm:staircase}
Let $S$ be a numerical set of genus $g\ge1$ with gaps $x_1<\dots<x_g$, let $\lambda=\lambda(S)$, and let $u=c_0c_1\cdots
c_F$ be its Keith--Nath word (so $c_0=1$ and $\cs(\lambda)=c_1\cdots c_F$ when the smallest part of $\lambda$ is $1$). The
following conditions are equivalent:
\begin{enumerate}
\item[(i)] $x_i\le2i-1$ for $i=1,\dots,g$;
\item[(ii)] $\lambda\subseteq\delta_g$;
\item[(iii)] for $i=1,\dots,g$, at least $i$ of the letters $c_1,\dots,c_{2i-1}$ are zeros, where
$c_x=1$ for $x>F$.
\end{enumerate}
If $S$ is a numerical semigroup, these conditions hold. Moreover, equality $x_i=2i-1$ for all $i$ (equivalently
$\lambda=\delta_g$) holds if and only if $S=\langle2,2g+1\rangle$.
\end{theorem}

\begin{proof}
Equivalence of (i) and (ii): by Lemma~\ref{lem:rowgap}, $\lambda_{g+1-i}\le g+1-(g+1-i)=i$ is equivalent to
$x_i-(i-1)\le i$, that is, to $x_i\le2i-1$.

Equivalence of (i) and (iii): the zeros of $u$ at positions $\ge1$ are exactly the gaps. So at least $i$ zeros among
$c_1,\dots,c_{2i-1}$ says that at least $i$ gaps lie in $[1,2i-1]$, that is, $x_i\le2i-1$.

Now, we show that every numerical semigroup satisfy (i). Let $S$ be a numerical semigroup and fix $i$. Partition $\{1,\dots,x_i-1\}$ into the
classes $\{y,x_i-y\}$; there are $\lceil(x_i-1)/2\rceil$ of them. If both elements of a class were in $S$, their sum
$x_i$ would be in $S$, which is false. So each class contains a gap smaller than $x_i$, and distinct classes are
disjoint. Since exactly $i-1$ gaps are smaller than $x_i$, we get $i-1\ge\lceil(x_i-1)/2\rceil\ge(x_i-1)/2$, that is,
$x_i\le2i-1$.

Finally, if $x_i=2i-1$ for all $i$, then the gaps are $1,3,\dots,2g-1$, so $S$ consists of all even numbers
and all integers $\ge2g$, that is, $S=\langle2,2g+1\rangle$. By Lemma~\ref{lem:rowgap}, this is the case
$\lambda_{g+1-i}=i$ for all $i$, that is, $\lambda=\delta_g$. Conversely, $\langle2,2g+1\rangle$ has gaps
$1,3,\dots,2g-1$.
\end{proof}

\begin{remark}
The case $i=g$ of Theorem~\ref{thm:staircase} is $F\le2g-1$, equivalently $\lambda_1\le\ell(\lambda)$. So
Theorem~\ref{thm:staircase} refines the first statement of Theorem~\ref{thm:ellbound}: it bounds every row, not only the
first one. The gap inequality $x_i\le2i-1$ is a classical argument for gaps. What is new here is its reformulation as
containment in a staircase and as a ballot condition on creation sequences, together with the consequences below. In
the threshold-graph language of Section~\ref{sec:cs}, the ballot condition says that for every $i\le g$, at least $i$
of the vertices $v_1,\dots,v_{2i-1}$ are isolated.
\end{remark}

\begin{corollary}\label{cor:weight}
Let $\lambda$ have smallest part $1$, and let $S=S_\lambda$ have genus $g$. Then
$
|\lambda|=g+w(S).
$
If $S$ is a numerical semigroup, then $0\le w(S)\le\binom g2$, and therefore
\[
g\le|\lambda|\le\frac{g(g+1)}{2}.
\]
Equality on the left holds exactly for $\lambda=(1^g)$ ($S=\{0\}\cup[g+1,\infty)$), and on the right exactly for
$\lambda=\delta_g$ ($S=\langle2,2g+1\rangle$). Moreover, the number of non-edges of the threshold graph $G(S)$ equals
$\sum_ix_i=g(g+1)/2+w(S)$.
\end{corollary}

\begin{proof}
By Lemma~\ref{lem:rowgap}, $|\lambda|=\sum_{i=1}^g(x_i-i+1)=\sum_ix_i-\binom g2=w(S)+\binom{g+1}2-\binom g2=w(S)+g$.
The gaps are distinct positive integers, so $x_i\ge i$ and $w(S)\ge0$, with equality if and only if the gaps are
$1,\dots,g$, that is, $\lambda=(1^g)$. For a semigroup, Theorem~\ref{thm:staircase} gives $x_i\le2i-1$, so
$w(S)\le\sum_i(i-1)=\binom g2$, with equality if and only if $x_i=2i-1$ for all $i$, that is, $\lambda=\delta_g$. Finally,
by Proposition~\ref{prop:graphdict}-(5), a pair $\{v_i,v_j\}$ with $i<j$ is a non-edge exactly when $j$ is
a gap, and each gap $x$ contributes $x$ such pairs.
\end{proof}

\begin{corollary}\label{cor:genusweight}
For all $n,g\ge1$, $\PG(n,g)$ is the number of numerical semigroups of genus $g$ and weight $n-g$. In particular:
\begin{enumerate}
\item $\PG(n,g)=0$ unless $g\le n\le g(g+1)/2$, and $\PG(g,g)=\PG\bigl(\tfrac{g(g+1)}2,g\bigr)=1$;
\item for $0\le w\le g$, the number of numerical semigroups of genus $g$ and weight $w$ is
\[
p(w)-E(2w-g),
\]
and it equals $p(w)$ when $2w\le g$.
\end{enumerate}
\end{corollary}

\begin{proof}
By Lemma~\ref{lem:bij}, the partitions of $n$ with $g$ parts corresponding to numerical semigroups are the
$\lambda(S)$ for the semigroups $S$ of genus $g$ with $|\lambda(S)|=n$, and $|\lambda(S)|=g+w(S)$ by
Corollary~\ref{cor:weight}. (1) follows from the bounds and equality cases in Corollary~\ref{cor:weight}. 

For (2), put $n=g+w$; the condition $w\le g$ is $2g\ge n$, and $2n-3g=2w-g$. Apply Theorem~\ref{thm:PG}, and
Lemma~\ref{lem:E}(1) when $2w\le g$.
\end{proof}

\begin{corollary}\label{cor:catalan}
Let $S$ be a numerical set of genus $g\ge1$ with Keith--Nath word $c_0c_1\cdots$ (with $c_x=1$ for $x>F$). Read the
$2g$ letters $c_1c_2\cdots c_{2g}$ as a lattice path, with $0$ an up-step and $1$ a down-step. Then $S$ satisfies the
conditions of Theorem~\ref{thm:staircase} if and only if this path is a Dyck path (it has $g$ up-steps and $g$
down-steps and never goes below its starting level). The map $S\mapsto c_1\cdots c_{2g}$ is a bijection from such
numerical sets onto Dyck paths of semilength $g$. Consequently there are exactly $C_g=\frac1{g+1}\binom{2g}g$ of them,
and the number $n_g$ of numerical semigroups of genus $g$ satisfies $n_g\le C_g$.
\end{corollary}

\begin{proof}
Put $a_i=x_i-i$. The up-steps of the path are at the positions of the gaps in $[1,2g]$. The $i$-th gap $x_i$ is preceded,
among $c_1,\dots,c_{x_i-1}$, by $i-1$ zeros and therefore by exactly $a_i=x_i-i$ ones.

Assume (i) of Theorem~\ref{thm:staircase}. Then all $g$ gaps lie in $[1,2g-1]$, so $c_1\cdots c_{2g}$ has exactly $g$
up-steps and $g$ down-steps. A path with $g$ up-steps and $g$ down-steps is a Dyck path if and only if for every
$i$, the $i$-th up-step is preceded by at most $i-1$ down-steps. (If some prefix had more down-steps than up-steps,
consider the first such prefix: it ends with a down-step, and the next up-step, say the $i$-th, is preceded by at least
$i$ down-steps; the converse is immediate.) Here this condition reads $a_i\le i-1$, which is (i). So the path is a Dyck
path.

Conversely, a Dyck path of semilength $g$ has $g$ up-steps, at positions $1\le x_1<\dots<x_g\le2g-1$ satisfying
$x_i-i\le i-1$. Putting $S=\N_0\setminus\{x_1,\dots,x_g\}$ gives the unique numerical set of genus $g$ satisfying (i)
whose word begins with the given path. Thus the map is a bijection. The number of Dyck paths of semilength $g$ is
$C_g$ \cite{Stanley}. The inequality $n_g\le C_g$ follows because every numerical semigroup satisfies (i)
of Theorem~\ref{thm:staircase}.
\end{proof}

\begin{remark}
For large $g$ the bound $n_g\le C_g$ is much weaker than the known bound $n_g\le1+3\cdot2^{g-3}$ ($g\ge3$)
\cite{BA}. The point of
Corollary~\ref{cor:catalan} is structural: the creation sequences of numerical semigroups of genus $g$ form a subset
of an explicit Catalan family defined by (iii) of Theorem 7.2. Numerically, $n_g/C_g=1,\,1,\,4/5,\,7/14,\,12/42,\dots$
for $g=1,\dots,5$.
\end{remark}

\section{The semigroup tree in partition language}\label{sec:tree}
The numerical semigroups form a rooted tree: the parent of $S\ne\N_0$ is $S\cup\{F(S)\}$, and the children of $S$ are
the sets $S\setminus\{x\}$ for the minimal generators $x>F(S)$ \cite{RGS}. Creation sequences turn this
tree into an operation on the largest part.

\begin{theorem}\label{thm:parent}
If $\lambda=(\lambda_1,\dots,\lambda_\ell)$ corresponds to a numerical semigroup and $\ell\ge2$, then so does
$(\lambda_2,\dots,\lambda_\ell)$, and it corresponds to $S\cup\{F(S)\}$. The partition $(1)$ corresponds to
$\{0,2,3,\dots\}$, whose parent is $\N_0\leftrightarrow\emptyset$.
\end{theorem}

\begin{proof}
Let $S=S^{KN}(\lambda)$ and $F=F(S)$. The set $S^\ast=S\cup\{F\}$ is a numerical semigroup: a sum of two nonzero
elements of $S^\ast$ involving $F$ exceeds $F$ and hence lies in $S$, and sums of elements of $S$ lie in $S$. Write
$u(\lambda)=u^\ast\,1^{\lambda_1-\lambda_2}\,0$, where $u^\ast$ ends with the north step of the second row. In
$S^\ast$ every integer beyond the last letter of $u^\ast$ is an element. So the Keith--Nath word of $S^\ast$ is
$u^\ast$, which is the word of $(\lambda_2,\dots,\lambda_\ell)$: the path of $\lambda$ with its top row removed.
\end{proof}

\begin{theorem}\label{thm:children}
Let $\lambda$ correspond to a numerical semigroup $S$ with multiplicity $m$ and Frobenius number $F$, and let
$d\ge0$. Then $(\lambda_1+d,\lambda_1,\lambda_2,\dots,\lambda_\ell)$ corresponds to a numerical semigroup if and only
if $F+d+1$ is a minimal generator of $S$. In particular this forces $d\le m-1$.
\end{theorem}

\begin{proof}
Prepending the part $\lambda_1+d$ appends $1^d0$ to the Keith--Nath word, so the new numerical set is
$S\setminus\{F+d+1\}$. For $x\in S\setminus\{0\}$, the set $S\setminus\{x\}$ is a numerical semigroup if and only if
$x$ is not a sum of two nonzero elements of $S\setminus\{x\}$. Since such summands are smaller than $x$, this is the
same as $x$ being a minimal generator. Every $y\ge F+m+1$ equals $m+(y-m)$ with $y-m>F$, so it is not a minimal
generator. Thus, $F+d+1\le F+m$.
\end{proof}

\begin{corollary}\label{cor:suffix}
If $\lambda$ corresponds to a numerical semigroup, then so does each of its suffixes
$(\lambda_i,\dots,\lambda_\ell)$. Every partition corresponding to a numerical semigroup is obtained from $(1)$ by
repeatedly prepending a new largest part $\lambda_1+d$, where $0\le d\le m-1$ and $F+d+1$ is a minimal generator.
\end{corollary}

\begin{proof}
Iterate Theorem~\ref{thm:parent}; the construction follows from Theorem~\ref{thm:children}.
\end{proof}

\section{Exact counts of semigroup partitions}\label{sec:count}

Here ``$\lambda$ corresponds to a numerical semigroup'' is meant in the sense of \cite{BNS}. By
Lemma~\ref{lem:bij}, it holds if and only if the smallest part of $\lambda$ is $1$ and $S_{\cs(\lambda)}$ is a
semigroup. Many of the partitions below have the shape
\begin{equation}\label{eq:shape}
\lambda=(t+1,\ \kappa_1+1,\ \dots,\ \kappa_r+1,\ 1^{m-1}),\qquad t\ge\kappa_1,\ m\ge2,
\end{equation}
where $\kappa=(\kappa_1,\dots,\kappa_r)$ is a partition. The parts of $\rho=(t+1,\kappa+1)$ are atleast 
$2$, so $\lambda$ has exactly $m-1$ parts equal to $1$ and $m(S^{KN}(\lambda))=m$. Its creation sequence is
\[
\cs(\lambda)=0^{m-1}\,w(\kappa)\,1^{t-\kappa_1}\,0 ,
\]
where $w(\kappa)$ occupies the positions $m,\dots,m+\kappa_1+r-1$, and $F=m+t+r$. Put $e=F-2m=t+r-m$.

\begin{lemma}\label{lem:key}
For $\lambda$ as in \eqref{eq:shape}:
\begin{enumerate}
\item if $e<0$, then $\lambda$ corresponds to a numerical semigroup;
\item if $e\ge0$ and $e\in\SSk$, then it does not;
\item if $e\ge0$, $e\notin\SSk$ and $m\ge\kappa_1+r$, then it does.
\end{enumerate}
\end{lemma}

\begin{proof}
Write $P$ for the set of $p\ge0$ such that $m+p\in S$ and $m+p<F$. By Lemma~\ref{lem:Sigma},
$P=(\Sigma(\kappa)\cap[0,\kappa_1+r-1])\cup[\kappa_1+r,\,t+r-1]$. The gaps larger than $m$ are $m+z_t$ for the zero
positions $z_t\le\kappa_1+r-1$ of $w(\kappa)$, together with $F=m+(t+r)$.

(1) If $e<0$ then $F<2m$, and Theorem~\ref{thm:short} applies.

(2) Since $m\ge2$ we have $e=t+r-m\le t+r-1$, so $P\cap[0,e]=\Sigma(\kappa)\cap[0,e]$. If $e=p+q$ with
$p,q\in\Sigma(\kappa)$, then $p,q\le e$ lie in $P$, and $(m+p)+(m+q)=2m+e=F$ is a gap.

(3) Suppose that $(m+p)+(m+q)$ is a gap for some $p,q\in P$. It is not $F$, since then $p+q=e$ with
$p,q\in P\cap[0,e]\subseteq\Sigma(\kappa)$. So $2m+p+q=m+z$ for a zero position $z\le\kappa_1+r-1<m$, giving
$m+p+q=z<m$, which is impossible. Sums involving $0$ or elements larger than $F$ are in $S$.
\end{proof}

\begin{remark}\label{rem:keygraph}
Let $G=G(S_\lambda)$ be the threshold graph of $\lambda$ (Definition~\ref{def:graphpartition}), with vertices
$v_0,\dots,v_F$. For $\lambda$ as in \eqref{eq:shape}, the construction of $G$ proceeds as follows:

(i) $v_1,\dots,v_{m-1}$ are isolated and $v_m$ is the first dominating vertex;

(ii) the vertices $v_m,\dots,v_{m+\kappa_1+r-1}$ are added according to the creation sequence $w(\kappa)$ of the
tail $\kappa$ ($v_{m+p}$ is dominating if and only if $p$ is a $1$-position of $w(\kappa)$);

(iii) a run of dominating vertices follows;

(iv) the last vertex $v_F=v_{2m+e}$ is isolated.

By Proposition~\ref{prop:graphdict}-(6), $S_\lambda$ fails to be a semigroup exactly when two dominating vertices
$v_x,v_y$ ($x,y\ge1$) have an isolated vertex $v_{x+y}$. Lemma~\ref{lem:key} says that, under its hypotheses, the only
isolated vertex that can be reached in this way is the last one, $v_F$. It is reached if and only if $e=p+q$ for two
$1$-positions $p,q$ of the tail's creation sequence (position $0$ and all positions beyond $w(\kappa)$ count as
$1$-positions). Thus whether $\lambda$ is exceptional is decided by the creation sequence of its tail alone.
This mechanism underlies all the solutions in Section~\ref{sec:solutions}.
\end{remark}

\subsection{Genus}

\begin{definition}\label{def:E}
For an integer $j$ let
\[
E(j)=\#\bigl\{\nu\ \text{a partition}:\ 2|\nu|+\ell(\nu)-j\notin\Sigma(\nu)\bigr\}.
\]
\end{definition}

\begin{lemma}\label{lem:E}
\begin{enumerate}
\item $E(j)=0$ for $j\le0$, and only partitions with $|\nu|<j$ contribute to $E(j)$; in particular $E(j)$ is finite.
\item $E(j)=\#\{\nu:\ 2|\nu|+\ell(\nu)<j\}+\#\{(\nu,t):\ 1\le t\le\ell(\nu),\ 2|\nu|-\nu_t+t=j\}$, and for each
$\nu$ at most one $t$ occurs.
\item Every $\nu$ counted by $E(j)$ satisfies $|\nu|+\ell(\nu)\le j$.
\end{enumerate}
\end{lemma}

\begin{proof}
Put $d=2|\nu|+\ell(\nu)-j$ and $r=\ell(\nu)$. By Lemma~\ref{lem:Sigma}, $d\notin\Sigma(\nu)$ if and only if $d<0$, or
$0\le d<\nu_1+r$ and $d=z_t=\nu_t+r-t$ for some $t$. The latter is the equation $2|\nu|-\nu_t+t=j$. The map
$t\mapsto2|\nu|-\nu_t+t$ is strictly increasing, so at most one $t$ occurs. This proves (2).

If $j\le0$ then $d\ge0$, and $d=z_t$ would give $2|\nu|-\nu_t+t\le0$, which is impossible for $\nu\ne\emptyset$
(since $|\nu|\ge\nu_t$ and $t\ge1$); for $\nu=\emptyset$ we have $\Sigma=\N_0\ni d$.
If $|\nu|\ge j$, then $d\ge|\nu|+r\ge\nu_1+r$, so $d\in\Sigma(\nu)$. This proves (1).

For (3): if $d<0$ then $|\nu|+r\le2|\nu|+r<j$. Otherwise $d=z_t<\nu_1+r$ gives $2|\nu|-\nu_1<j$, and
$|\nu|\ge\nu_1+r-1$ gives $|\nu|+r\le2|\nu|-\nu_1+1\le j$.
\end{proof}

\begin{theorem}\label{thm:PG}
For all integers $n\ge1$ and $g$ with $n/2\le g\le n$,
\[
\PG(n,g)=p(n-g)-E(2n-3g).
\]
\end{theorem}

\begin{proof}
Let $k=n-g$, so $k\le g$, and put $j=2k-g=2n-3g$.

Removing the first column is a bijection from partitions $\lambda\vdash n$ with exactly $g$
parts onto partitions $\mu\vdash k$ with $\ell(\mu)\le g$; the latter condition is automatic because
$\ell(\mu)\le k\le g$. So there are $p(k)$ candidates. Let $A=\ell(\mu)$, let $\nu$ be $\mu$ with its first column
removed, and put $s=|\nu|=k-A$ and $r=\ell(\nu)\le A$. Then
\[
\lambda=(\nu_1+2,\dots,\nu_r+2,\ 2^{A-r},\ 1^{g-A}),
\]
and $\lambda$ is determined by $\nu$ (given $n$ and $g$). A partition $\nu$ arises in this way if and only if
$A=k-s\ge r$, that is, $|\nu|+\ell(\nu)\le k$. We consider the following cases:

\textbf{Case-I: $A=g$}. Then $k=g$, $\mu=(1^g)$, $\nu=\emptyset$ and $\lambda=(2^g)$. This has no part $1$, so it does
not correspond to a semigroup (Lemma~\ref{lem:bij}), and it is counted by $E(j)=E(g)$ because
$-g\notin\Sigma(\emptyset)$ and $g\ge1$.

\textbf{Case-II: $A<g$}. Then $\lambda$ has $g-A\ge1$ parts equal to $1$, and
\[
\cs(\lambda)=0^{g-A}\,1\,0^{A-r}\,w(\nu).
\]
So $m=g-A+1$ is the multiplicity, the positions $m+1,\dots,T-1$ are gaps with
\[
T=m+(A-r)+1=g-r+2,
\]
and $w(\nu)$ occupies the positions $T,\dots,F$, where $F=T+\nu_1+r-1=g+\nu_1+1$. (If $\nu=\emptyset$, then $F=T-1$.)
The nonzero elements of $S$ up to $F$ are $m$ and $T+p$ for $p\in\Sigma(\nu)\cap[0,\nu_1+r-1]$. We check which
sums of two such elements can be gaps.\\
 
\textit{Two elements $\ge T$}: their sum is at least $2T=2g-2r+4$. Since $g\ge k=s+A\ge s+r$ and
$s\ge\nu_1+r-1$, we get $g\ge\nu_1+2r-1$, and so $2T\ge g+\nu_1+3>F$.\\

\textit{The element $m$ and an element $\ge T$}: the sum is at least
$m+T=2g-k+s-r+3\ge g+s-r+3$ (using $k\le g$). Since $s-r\ge\nu_1-1$, this exceeds $F=g+\nu_1+1$.\\

\textit{The element $m$ twice}: $2m=T+d$ with
\[
d=2m-T=2(g-k+s+1)-(g-r+2)=2s+r-(2k-g)=2|\nu|+\ell(\nu)-j .
\]
Thus $S$ is a semigroup if and only if $2m$ is not a gap. If $d<0$, then $m<2m<T$ and $2m$ is a gap. If $d\ge0$, then
$2m$ is a gap if and only if $d$ is the position of a zero of $w(\nu)$. The positions $\ge\nu_1+r$ correspond to
integers $>F$, which lie in $S$. So in all cases $\lambda$ fails if and only if $d\notin\Sigma(\nu)$.

The failing $\lambda$ correspond bijectively to the partitions $\nu$ with $d\notin\Sigma(\nu)$ that
arise, specifically with $|\nu|+\ell(\nu)\le k$. By Lemma~\ref{lem:E}-(3) every $\nu$ counted by $E(j)$ satisfies
$|\nu|+\ell(\nu)\le j\le k$, so all of them arise. Therefore exactly $E(j)$ candidates fail.
\end{proof}

\begin{proposition}\label{prop:GF}
The generating function of $E$ is
\[
\sum_{j\ge1}E(j)x^j=\frac{x}{1-x}\prod_{i\ge1}\frac1{1-x^{2i+1}}
+\sum_{v\ge1}\sum_{t\ge1}\frac{x^{(2t-1)v+t}}{(x^2;x^2)_{t-1}\,(x^2;x^2)_{v}},
\]
where $(q;q)_s=\prod_{i=1}^{s}(1-q^i)$.
\end{proposition}

\begin{proof}
We use Lemma~\ref{lem:E}-(2). Since $2|\nu|+\ell(\nu)=\sum_i(2\nu_i+1)$, the generating function of
$2|\nu|+\ell(\nu)$ over all $\nu$ is $\prod_{i\ge1}(1-x^{2i+1})^{-1}$. Multiplying by $x/(1-x)=\sum_{j\ge1}x^j$
counts, for each $j$, the partitions with $2|\nu|+\ell(\nu)<j$. For the second term, fix $t\ge1$ and $v=\nu_t\ge1$.
The parts above row $t$ are $v+\alpha_1\ge\dots\ge v+\alpha_{t-1}$ with $\alpha$ a partition into at most $t-1$
parts. The parts below row $t$ form a partition $\beta$ with all parts $\le v$. Then
$$2|\nu|-\nu_t+t=2(t-1)v+2|\alpha|+v+2|\beta|+t=(2t-1)v+t+2|\alpha|+2|\beta|.$$ The generating functions of $\alpha$
and $\beta$, weighted by $x^{2|\cdot|}$, are $1/(x^2;x^2)_{t-1}$ and $1/(x^2;x^2)_v$.
\end{proof}

As the consequence of the preceding result we recover Theorems 1, 2, and 5 \cite{BNS},   and provide the correct versions of Theorem 3 and Corollary 4. 

\begin{corollary}\label{cor:PG}
\begin{enumerate}
\item $\PG(n,g)=p(n-g)$ whenever $3g\ge2n$ \cite[Theorem 1]{BNS}. The bound is sharp: if $n/2\le g<2n/3$, then
$\PG(n,g)<p(n-g)$.
\item $\PG(3n-i,2n-i)=p(n)-E(i)$ for $n\ge i$. Since $E(1)=1$ and $E(2)=E(3)=2$, this gives \cite[Theorem 2]{BNS}.
\item $\PG(3n+2-k,2n-k)=p(n+2)-4-k=b(n)-k$ for $k=0,1,2$ and $n\ge k+2$. This corrects \cite[Theorem 3, Corollary 4]{BNS},
which state $b(n)-(k+2)$ and $p(n+2)-6-k$.
\item $\PG(3n+2-j,2n-j)-\PG(3n-1-j,2n-2-j)=\bar p(n+2)$ for $n\ge j+3$, which is \cite[Theorem 5]{BNS} with the
improved range $n\ge j+3$ in place of $n\ge j+5$.
\end{enumerate}
\end{corollary}

\begin{proof}
(1) If $3g\ge2n$ then $2n-3g\le0$ and $E=0$ by Lemma~\ref{lem:E}-(1). If $j=2n-3g\ge1$, then $\nu=\emptyset$ is counted by $E(j)$ (since $-j<0$), so $E(j)\ge1$.

(2) For $(3n-i,2n-i)$ we have $2N-3g=i$, and $g\ge N/2$ if and only if $n\ge i$. The values of $E$ come from
Table~\ref{tab:EH}. By hand, with $d=2|\nu|+\ell(\nu)-j$: $E(1)$ counts only $\emptyset$ (for $\nu=(1)$, $d=2$ is
beyond $w((1))=10$). $E(2)$ counts $\emptyset$ and $(1)$ ($d=1$ is the zero of $10$); for $(2)$ and $(1,1)$ we get
$d=3,4$, beyond words of length $3$. $E(3)$ counts $\emptyset$ and $(2)$ ($d=2$ is the zero of $w((2))=110$), while
$(1)$ gives $d=0\in\Sigma$ and $(1,1)$ gives $d=3$.

(3) Here $2N-3g=k+4$, and $g\ge N/2$ if and only if $n\ge k+2$. From Table~\ref{tab:EH}, $E(4)=4$, $E(5)=5$ and
$E(6)=6$, so $\PG=p(n+2)-4-k$. Moreover $b(n)=p(n+2)-4$: among the partitions of $n+2$, exactly
$(1^{n+2})$ and $(2,1^n)$ have more than $n$ parts, and exactly $(n+2)$ and $(n+1,1)$ have a part larger than $n$.
In the proof of \cite[Theorem 3]{BNS}, the partitions $\{n+2,2,2,1^{2n-k-2}\}$ and $\{n+3,1^{2n-k-1}\}$ are removed
from the domain of the map but not added back, although for $n\ge7$ they correspond to numerical semigroups (for
instance, the second is a hook with $d=n+2\le2n-k-1$, which is covered by Corollary~\ref{cor:hook}).

(4) For both terms $2N-3g=j+4$, so the $E$-terms cancel, and $p(n+2)-p(n+1)=\bar p(n+2)$. The conditions
$g\ge N/2$ for the two terms are $n\ge j+2$ and $n\ge j+3$.
\end{proof}

\begin{corollary}\label{cor:rem1}
The exceptional-case functions of \cite[Remark 1]{BNS} are given by
\begin{align*}
h_1(i)&=p(n)-\PG(3n-i,2n-i)=E(i) && (n\ge i),\\
h_2(i)&=b(n)-\PG(3n+2-i,2n-i)=E(i+4)-4 && (n\ge i+2),
\end{align*}
and they are independent of $n$. More generally, for any line $\PG(L_1(n)-k,L_2(n)-k)$ with
$2L_2(n)-k\ge L_1(n)$, we have
\[
p(L_1(n)-L_2(n))-\PG(L_1(n)-k,L_2(n)-k)=E(2L_1(n)-3L_2(n)+k).
\]
\end{corollary}

\begin{proof}
All three identities are Theorem~\ref{thm:PG} with the stated arguments. For the second, use
$b(n)=p(n+2)-4$ from Corollary~\ref{cor:PG}-(3).
\end{proof}

\subsection{Multiplicity}

\begin{definition}\label{def:H}
For $i\ge0$, let
\[
H(i)=\#\bigl\{\kappa\ \text{a partition}:\ |\kappa|\le i,\ \ i-|\kappa|\in\SSk\bigr\},
\]
and set $H(i)=0$ for $i<0$. The empty partition contributes to $H(i)$ for every $i\ge0$.
\end{definition}

\begin{theorem}\label{thm:PM}
For all integers $m\ge2$ and $n\ge m-1$ with $3m>n$,
\[
\PM(n,m)=\bar p(n-m+1)-H(n-2m).
\]
\end{theorem}

\begin{proof}
If $n=m-1$, the only candidate is $\lambda=(1^{m-1})$, which corresponds to $\langle m,m+1,\dots,2m-1\rangle$, and
$\bar p(0)=1$; if $n=m$, no such $\lambda$ exists and $\bar p(1)=0$. So assume $n\ge m+1$.
Put $i=n-2m$ and $R=n-m+1$.

A partition of $n$ with exactly $m-1$ parts equal to $1$ is $\lambda=(\rho,1^{m-1})$ with
$\rho\vdash R$ having no part $1$; there are $\bar p(R)$ of them. Writing
$\rho=(t+1,\kappa_1+1,\dots,\kappa_r+1)$ gives a bijection between such $\rho$ and pairs $(t,\kappa)$ with
$t\ge\kappa_1$ and $t+1+|\kappa|+r=R$. Each $\lambda$ has the shape \eqref{eq:shape}, and
\[
e=t+r-m=(R-1-|\kappa|-r)+r-m=i-|\kappa| .
\]
Given $\kappa$, the value of $t$ is determined, and $\kappa$ arises if and only if $t=m+e-r\ge\kappa_1$.

If $|\kappa|>i$, then $e<0$ and $\lambda$ corresponds to a semigroup by
Lemma~\ref{lem:key}-(1). Let $|\kappa|\le i$. Since $3m\ge n+1$ we have $i\le m-1$, so
$\kappa_1+r\le|\kappa|+1\le m$ (using $|\kappa|\ge\kappa_1+r-1$). Therefore  Lemma~\ref{lem:key}-(2),(3) apply, and
$\lambda$ fails if and only if $i-|\kappa|\in\SSk$. Moreover every such $\kappa$ arises, since
$t=m+e-r\ge m-r\ge\kappa_1$.

Thus the failing $\lambda$ are in bijection with the partitions counted by $H(i)$.
\end{proof}

\begin{theorem}\label{thm:PMboundary}
For every $m\ge2$,
\[
\PM(3m,m)=\bar p(2m+1)-H(m)+m .
\]
In particular Theorem~\ref{thm:PM} fails at $n=3m$, and the range $3m>n$ is sharp.
\end{theorem}

\begin{proof}
Now $i=m$ and $R=2m+1$, and we reuse the parametrisation from the proof of Theorem~\ref{thm:PM}. Partitions $\kappa$
with $|\kappa|>m$ give semigroups, as before. Let $|\kappa|\le m$, so $e=m-|\kappa|\ge0$. Then $\kappa$ arises if and
only if $t=2m-|\kappa|-r\ge\kappa_1$, that is, $\kappa_1+r+|\kappa|\le2m$. Since $\kappa_1+r\le|\kappa|+1$, with equality
exactly for hooks, this fails precisely when $|\kappa|=m$ and $\kappa$ is a hook. There are exactly $m$ hooks of size
$m$. For every other $\kappa$ with $|\kappa|\le m$ we have $\kappa_1+r\le m$, so Lemma~\ref{lem:key}-(2),(3) apply
exactly as before.

Therefore the failing $\lambda$ correspond to the $\kappa$ counted by $H(m)$, except the $m$ hooks of size $m$, which do
not arise. Each of these hooks is counted by $H(m)$, because $e=0\in\SSk$. Therefore the number of failures is
$H(m)-m$.
\end{proof}

\begin{theorem}\label{thm:PMboundary2}
For every $m\ge2$,
\[
\PM(3m+1,m)=\bar p(2m+2)-H(m+1)+m .
\]
\end{theorem}

\begin{proof}
Now $i=m+1$ and $R=2m+2$, and we again use the parametrisation $\rho\leftrightarrow(t,\kappa)$ from the proof of
Theorem~\ref{thm:PM}, with $e=m+1-|\kappa|$ and $t=m+e-r$. Partitions $\kappa$ with $|\kappa|>m+1$ have $e<0$ and
give semigroups. Let $|\kappa|\le m+1$. Then $\kappa$ arises if and only if $t\ge\kappa_1$, that is,
\[
\kappa_1+r+|\kappa|\le2m+1 .
\]
Recall that $\kappa_1+r\le|\kappa|+1$, with equality exactly for hooks, and that $|\kappa|-(\kappa_1+r-1)$ is the
number of cells of $\kappa$ outside its first row and first column. We sort the partitions with $|\kappa|\le m+1$ into
four classes.

\emph{(a) Not arising.} The displayed inequality fails exactly when $|\kappa|=m+1$ and $\kappa_1+r\ge m+1$. These are
the $m+1$ hooks of size $m+1$, and the $m-2$ partitions $(a,2,1^{b})$ with $a\ge2$ and $a+b=m-1$ (one cell outside the
first row and column). For all of them $e=0\in\SSk$, so all $2m-1$ are counted by $H(m+1)$ although they do not
arise.

\emph{(b) Arising with $\kappa_1+r\le m$.} Here $m\ge\kappa_1+r$, and Lemma~\ref{lem:key}-(2),(3) show that $\lambda$
fails if and only if $e\in\SSk$, exactly as counted by $H(m+1)$.

\emph{(c) Arising with $\kappa_1+r\ge m+1$.} Since $|\kappa|\le m+1$, and the case $|\kappa|=m+1$ belongs to (a), this
forces $|\kappa|=m$ and $\kappa$ a hook, $\kappa=(m-s,1^{s})$ with $0\le s\le m-1$, and then $e=1$. For $s=0$, $1\in
\Sigma((m))$, so $\lambda$ fails by Lemma~\ref{lem:key}-(2), and it is counted by $H(m+1)$. For $s\ge1$,
$\Sigma(\kappa)=\{0\}\cup[s+1,m-1]\cup[m+1,\infty)$ (see the proof of Lemma~\ref{lem:hooksums}), so
$1\notin\SSk$ and $H(m+1)$ does not count $\kappa$. Nevertheless $\lambda$ fails: the last letter of $w(\kappa)$ is a
$0$ at position $\kappa_1+r-1=m$, so $2m=m+m$ is a gap while $m\in S$. This gives $m-1$ failures not counted by
$H(m+1)$.

The number of failures is therefore $H(m+1)-(2m-1)+(m-1)=H(m+1)-m$, and
$\PM(3m+1,m)=\bar p(2m+2)-H(m+1)+m$.
\end{proof}

In the following result, we recover and obtain the correct versions of some results discussed in \cite{BNS}, Section 5.

\begin{corollary}\label{cor:PM}
Put $h(i)=H(i)-1$. Since $a(n+1)=\bar p(n+1)-1$, we have $a(n+1)-\PM(2n-i,n-i)=h(i)$ whenever $n\ge2i+1$.
\begin{enumerate}
\item If $m>n/2$ (with $m\ge2$ and $n\ge m-1$), then $\PM(n,m)=\bar p(n+1-m)$ \cite[Proposition 9]{BNS}.
\item $\PM(2n-i,n-i)=a(n+1)-h(i)$ for $n\ge2i+1$, with $h(0),h(1),h(2),h(3)=0,1,2,5$. This is
\cite[Theorem 10]{BNS}. For $i=3$ the correct range is $n\ge7$: at $n=6$ we have $3m=n$, and
Theorem~\ref{thm:PMboundary} gives $\PM(9,3)=a(7)-5+3=a(7)-2$. At $n=5$, direct computation gives
$\PM(7,2)=0=a(6)-3$.
\item $\PM(2n-j,n-j)-\PM(2(n-1)-j,(n-1)-j)=\bar p(n+1)-\bar p(n)=\bar a(n+1)$ for $n\ge2j+2$ \cite[Theorem 12,
Corollary 13]{BNS}.
\item $\PM(3j+2,j+1)-\PM(3j,j)=\bar a(2j+2)-j$ for $j\ge2$ \cite[Theorem 14]{BNS}.
\item $\PM(3j,j)-\PM(3j-2,j-1)=\bar a(2j+1)+1$ for $j\ge3$ \cite[Theorem 15]{BNS}.
\end{enumerate}
\end{corollary}

\begin{proof}
(1) If $m>n/2$, then $i<0$, so $H(i)=0$, and $3m>n$.

(2) For $(2n-i,n-i)$ the condition $3m>n$ reads $n>2i$. The values $H(0),\dots,H(3)=1,2,3,6$ are computed directly
from Definition~\ref{def:H}; for example, $H(1)$ counts $\emptyset$ (as $1\in\N_0$) and $(1)$ (as $0\in\SSk$). For
$n=5$: a partition of $7$ with exactly one part $1$ has a unique smallest part, so by Theorem~\ref{thm:unique} it
corresponds to a semigroup only if it is a staircase $(k,\dots,2,1)$ of size $7$; there is none, so
$\PM(7,2)=0=a(6)-3$.

(3) Both terms have $i=j$, so the $H$-terms cancel. For $n\ge2$ the partitions of $n+1$ with no part $1$ split into
those with $\ell\ge2$ and $\lambda_1=\lambda_2$, counted by $\bar a(n+1)$, and those with $\lambda_1>\lambda_2$ (or
$\ell=1$). Subtracting $1$ from the largest part maps the latter bijectively onto the partitions of $n$ with no part
$1$; the inverse adds $1$ to the largest part. Thus $\bar p(n+1)-\bar p(n)=\bar a(n+1)$ (compare
\cite[Corollary 13]{BNS}). The condition $3m>n$ for the second term reads $n\ge2j+2$.

(4) Since $3(j+1)>3j+2$, Theorem~\ref{thm:PM} gives $\PM(3j+2,j+1)=\bar p(2j+2)-H(j)$, and
Theorem~\ref{thm:PMboundary} gives $\PM(3j,j)=\bar p(2j+1)-H(j)+j$. Subtract, and use (3).

(5) Theorem~\ref{thm:PMboundary} gives $\PM(3j,j)=\bar p(2j+1)-H(j)+j$, and Theorem~\ref{thm:PMboundary2} with
$m=j-1$ gives $\PM(3j-2,j-1)=\bar p(2j)-H(j)+j-1$. Subtract, and use (3).
\end{proof}

The following presents an answer to \cite{BNS}, Remark 4.

\begin{corollary}\label{cor:rem4}
The function asked for in \cite[Remark 4]{BNS} is $h(i)=H(i)-1$, valid for $n\ge2i+1$. Its first values are
listed in Table~\ref{tab:EH}.
\end{corollary}

\begin{proof}
This is Corollary~\ref{cor:PM}-(2), Theorem~\ref{thm:PM} applied to $(2n-i,n-i)$ gives
$\PM(2n-i,n-i)=\bar p(n+1)-H(i)=a(n+1)-(H(i)-1)$ whenever $3(n-i)>2n-i$, that is, $n\ge2i+1$.
\end{proof}

\subsection{Frobenius number}

For a partition $\pi$ put $c_\pi=\pi_1+\ell(\pi)$ (so $c_\emptyset=0$), and for $\varepsilon\in\{0,1\}$ let
\[
G_\pi(\varepsilon)=\#\bigl\{e:\ 0\le e<c_\pi,\ e\equiv\varepsilon\ (\mathrm{mod}\ 2),\
e\notin\Sigma(\pi)+\Sigma(\pi)\bigr\}.
\]

\begin{theorem}\label{thm:PF}
For all integers $k\ge0$ and $n\ge4k+1$,
\[
\PF(n,n-k)=\sum_{\pi\vdash k}\Bigl(\Bigl\lceil\frac{n-k}{2}\Bigr\rceil-c_\pi+G_\pi(n-k)\Bigr)
=p(k)\Bigl\lceil\frac{n-k}{2}\Bigr\rceil-C_k(n-k),
\]
where $C_k(\varepsilon)=\sum_{\pi\vdash k}\bigl(c_\pi-G_\pi(\varepsilon)\bigr)$ depends only on $k$ and on the parity of
$n-k$ (Table~\ref{tab:PF}).
\end{theorem}

\begin{proof}
Let $\lambda\vdash n$ have $h_{(1,1)}(\lambda)=n-k$. Let $t=\lambda_1-1$ be its arm and
$b=\ell(\lambda)-1$ its leg, so $t+b+1=n-k$. The cells outside the first row and the first column form the partition
$\pi=(\lambda_2-1,\dots,\lambda_\ell-1)$ (zeros omitted), with $|\pi|=k$, $r=\ell(\pi)\le b$ and $\pi_1\le t$.
Conversely every triple $(t,b,\pi)$ with $\pi\vdash k$, $t\ge\pi_1$, $b\ge r$ and $t+b=n-k-1$ arises from exactly one
$\lambda=(t+1,\pi_1+1,\dots,\pi_r+1,1^{b-r})$. This $\lambda$ can correspond to a semigroup only if it has a part $1$,
i.e.\ $b\ge r+1$. (For $k=0$ and $n=1$ the only partition is $\lambda=(1)$, which corresponds to
$\mathbb{N}_0\setminus\{1\}$, in agreement with $\lceil 1/2\rceil=1$; from now on we assume $n-k\ge2$.) In that case it has the shape \eqref{eq:shape} with $\kappa=\pi$ and $m=b-r+1\ge2$, and
\[
e=t+r-m=n-k-2-2b+2r,\qquad m=\frac{n-k-e}{2}.
\]
For fixed $\pi$, as $b$ runs from $r+1$ to $n-k-1-\pi_1$ (the upper limit is $t\ge\pi_1$),
$e$ runs through all integers $e\equiv n-k\pmod 2$ in the interval
\[
I_\pi=\bigl[-(n-k)+2c_\pi,\ n-k-4\bigr].
\]

By Lemma~\ref{lem:key}, the classification of $e\in I_\pi$ is given as follows.

(i) $e<0$: semigroup. The number of such $e$ of the right parity in $I_\pi$ is
$\lceil(n-k)/2\rceil-c_\pi$. Indeed, they are $-(n-k)+2c_\pi+2s$ for $0\le s\le(n-k-2c_\pi-1)/2$, and
$\lfloor(n-k-1)/2\rfloor+1=\lceil(n-k)/2\rceil$. This needs $n-k\ge2c_\pi$, which holds because $c_\pi\le k+1$ and
$n\ge4k+1$ (for $k=0$, $c_\pi=0$).

(ii) $e\ge c_\pi$: failure, by Lemma~\ref{lem:key}(2), because $e\in\Sigma(\pi)\subseteq\Sigma(\pi)+\Sigma(\pi)$
(Lemma~\ref{lem:Sigma}).

(iii) $0\le e<c_\pi$: failure if $e\in\Sigma(\pi)+\Sigma(\pi)$, and a semigroup otherwise, provided
$m\ge c_\pi$ (Lemma~\ref{lem:key}(3)). For these $e$ we have $m=(n-k-e)/2\ge\lceil(n-k-c_\pi+1)/2\rceil$, and this is
at least $c_\pi$ as soon as $n\ge k+3c_\pi-2$. That holds since $c_\pi\le k+1$ and $n \ge 4k+1$.

It remains to check that all $e\in[0,c_\pi)$ of the right parity lie in $I_\pi$. The lower end is fine since
$-(n-k)+2c_\pi\le0$. For the upper end, $c_\pi-1\le k\le n-k-4$ when $k\ge2$. For $k=1$ (so $\pi=(1)$ and $c_\pi=2$)
the right parity and $n\ge5$ give $e\le n-5$ in all cases. For $k=0$ there is nothing to check.

For each $\pi$ we obtain $\lceil(n-k)/2\rceil-c_\pi+G_\pi(n-k)$ semigroups; now sum over
$\pi\vdash k$.
\end{proof}

\begin{lemma}\label{lem:hooksums}
Let $k\ge2$ and let $\pi=(k-j,1^j)$ be a hook of size $k$, $0\le j\le k-1$. Then $k\notin\Sigma(\pi)+\Sigma(\pi)$ if
and only if $j\ge\lceil(k-1)/2\rceil$. Consequently exactly $\lceil k/2\rceil$ hooks of size $k$ satisfy
$k\notin\Sigma(\pi)+\Sigma(\pi)$.
\end{lemma}

\begin{proof}
For $j=0$, $w(\pi)=1^k0$ and $\Sigma(\pi)=[0,k-1]\cup[k+1,\infty)$, so $k=1+(k-1)$ lies in the sumset. For $j\ge1$,
$w(\pi)=1\,0^{j}\,1^{k-j-1}\,0$, so $\Sigma(\pi)=\{0\}\cup[j+1,k-1]\cup[k+1,\infty)$. Writing $k=x+y$ with
$x,y\in\Sigma(\pi)$: if $x=0$, then $y=k\notin\Sigma(\pi)$; otherwise $x,y\in[j+1,k-1]$, which is possible if and
only if $2(j+1)\le k$. So $k\notin\Sigma(\pi)+\Sigma(\pi)$ exactly when $2j+2>k$, that is,  $j\ge\lceil(k-1)/2\rceil$. The number of such
$j\in[0,k-1]$ is $k-\lceil(k-1)/2\rceil=\lceil k/2\rceil$.
\end{proof}

\begin{theorem}\label{thm:PFboundary}
For every $k\ge1$,
\[
\PF(4k,3k)=p(k)\Bigl\lceil\frac{3k}{2}\Bigr\rceil-C_k(k\bmod2)-\Bigl\lceil\frac k2\Bigr\rceil .
\]
In particular, Theorem~\ref{thm:PF} fails at $n=4k$, and the range $n\ge4k+1$ is sharp.
\end{theorem}

\begin{proof}
For $k=1$: the partitions of $4$ with largest hook $3$ are just $(2,2)$, which has no part $1$, so $\PF(4,3)=0$. The
right-hand side is $1\cdot2-1-1=0$.

Let $k\ge2$ and $n=4k$, and follow the proof of Theorem~\ref{thm:PF}. The conditions $n-k\ge2c_\pi$ and
$c_\pi-1\le n-k-4$ still hold, since $3k\ge2k+2$ and $3k-4\ge k$. The only step that can fail is $m\ge c_\pi$ for
$0\le e<c_\pi$. Here $m=(3k-e)/2$, so $m<c_\pi$ implies $e>3k-2c_\pi$. Since $e<c_\pi$, this requires
$3k-2c_\pi<c_\pi-1$, which implies $c_\pi>k+\tfrac13$, that is, $c_\pi=k+1$. This implies $\pi$ is a hook of size $k$. For a hook,
$e\equiv3k\equiv k\pmod2$, and $3k-2(k+1)<e\le k$ forces $e=k$ and $m=k$.

So the formula of Theorem~\ref{thm:PF} is correct except possibly for the pairs ($\pi$ a hook, $e=k$). If
$k\in\Sigma(\pi)+\Sigma(\pi)$, then the formula already counts this pair as a failure by
Lemma~\ref{lem:key}-(2). If $k\notin\Sigma(\pi)+\Sigma(\pi)$, the formula counts it as a semigroup. However $m=k$,
and position $k=\pi_1+r-1$ is the last letter of $w(\pi)$, which is $0$. So $2m=m+k$ is a gap while $m\in S$, and
$\lambda$ is not a semigroup. By Lemma~\ref{lem:hooksums} there are exactly $\lceil k/2\rceil$ such hooks, and the
formula overcounts by exactly this number.
\end{proof}

In the following result we recover and correct some results discussed in \cite{BNS}, Section 4.

\begin{corollary}\label{cor:PF}
\begin{enumerate}
\item $\PF(n,n)=\lceil n/2\rceil$ \cite[Theorem 6]{BNS}.
\item $\PF(2n,2n-1)=n-1$ and $\PF(2n+1,2n)=n-2$ for $n\ge3$ \cite[Theorem 7-(1),(2)]{BNS}.
\item $\PF(2n,2n-2)=\PF(2n-1,2n-3)=2n-7$ for $n\ge5$; and $\PF(2n,2n-3)=3n-11$, $\PF(2n-1,2n-4)=3n-17$ for $n\ge7$.
These are \cite[Theorem 7(3)--(5)]{BNS}, whose hypothesis $n>2$ must be strengthened as stated. For example
$\PF(6,4)=0\ne-1$, and $\PF(12,9)=5\ne7$ by Theorem~\ref{thm:PFboundary}.
\item $\PF(n,n-k)-\PF(n-2,n-k-2)=p(k)$ for $n\ge4k+3$ \cite[Theorem 8]{BNS}.
\end{enumerate}
\end{corollary}

\begin{proof}
We use the constants $C_0=0$, $C_1(0)=2$, $C_1(1)=1$, $C_2(0)=C_2(1)=5$, $C_3(0)=11$ and $C_3(1)=8$. We derive two
of them. 

For $k=1$, $\pi=(1)$, $c_\pi=2$ and $\Sigma(\pi)=\{0\}\cup[2,\infty)$, so $0$
lies in the sumset and $1$ does not; therefore, $G(0)=0$, $G(1)=1$, $C_1(0)=2$ and $C_1(1)=1$. The others follow in the
same way from Definition~\ref{def:KNset}.

(1) is Corollary~\ref{cor:hook}.

(2) Take $k=1$. For $\PF(2n,2n-1)$ the quantity $n'-k=2n-1$ is odd, so $\PF=\lceil(2n-1)/2\rceil-C_1(1)=n-1$. For
$\PF(2n+1,2n)$ it is $2n$, which is even, so $\PF=n-C_1(0)=n-2$. The range $n'\ge4k+1=5$ holds for $n\ge3$.

(3) $k=2$: $\PF=2\lceil(n'-2)/2\rceil-5$, giving $2(n-1)-5=2n-7$ for $n'=2n$ and for $n'=2n-1$; the range $n'\ge9$
means $n\ge5$. For $k=3$: $\PF(2n,2n-3)=3(n-1)-C_3(1)=3n-11$ and $\PF(2n-1,2n-4)=3(n-2)-C_3(0)=3n-17$; the range
$n'\ge13$ means $n\ge7$. At the excluded values: the only partitions of $6$ with largest hook $4$ are $(3,3)$ and
$(2,2,2)$, neither of which has a part $1$, so $\PF(6,4)=0$; and
$\PF(12,9)=3\cdot5-8-2=5$ by Theorem~\ref{thm:PFboundary}.

(4) When $n$ increases by $2$, $\lceil(n-k)/2\rceil$ increases by $1$ and the parity is unchanged. The range
$n-2\ge4k+1$ is $n\ge4k+3$.
\end{proof}

In the following result we provide an answer to \cite{BNS}, Remark 3.

\begin{corollary}\label{cor:rem3}
For $2n\ge4j+1$,
\[
\PF(2n,2n-j)=p(j)\Bigl\lceil n-\frac j2\Bigr\rceil-C_j(j\bmod 2),
\]
which replaces the upper bound of \cite[Remark 3]{BNS} by an identity.
\end{corollary}

\begin{proof}
Apply Theorem~\ref{thm:PF} with $2n$ in place of $n$ and $k=j$. Then $\lceil(2n-j)/2\rceil=\lceil n-j/2\rceil$, and
$2n-j\equiv j\pmod2$.
\end{proof}

\subsection{Tables}

\begin{table}[ht]
\centering
\small
\begin{tabular}{c|rrrrrrrrrrrrrrr}
$j$ & 1&2&3&4&5&6&7&8&9&10&11&12&13&14&15\\ \hline
$E(j)$ & 1&2&2&4&5&6&10&12&15&20&26&32&40&50&60\\ \hline\hline
$j$ & 16&17&18&19&20&21&22&23&24&25&26&27&28&29&30\\ \hline
$E(j)$ & 76&92&110&134&160&191&230&272&320&380&446&522&612&715&830
\end{tabular}

\vspace{8pt}
\begin{tabular}{c|rrrrrrrrrrrrr}
$i$ & 0&1&2&3&4&5&6&7&8&9&10&11&12\\ \hline
$H(i)$ & 1&2&3&6&9&15&21&33&46&69&91&132&176\\ \hline
$h(i)$ & 0&1&2&5&8&14&20&32&45&68&90&131&175
\end{tabular}
\caption{The exceptional-case functions $E$ (genus) and $H$, $h=H-1$ (multiplicity).}
\label{tab:EH}
\end{table}

\begin{table}[ht]
\centering
\small
\begin{tabular}{c|rrrrrrrrrrr}
$k$ & 0&1&2&3&4&5&6&7&8&9&10\\ \hline
$p(k)$ & 1&1&2&3&5&7&11&15&22&30&42\\
$C_k(0)$ & 0&2&5&11&20&36&58&95&145&221&323\\
$C_k(1)$ & 0&1&5&8&19&29&55&81&134&194&301
\end{tabular}
\caption{Constants in Theorem~\ref{thm:PF}: $\PF(n,n-k)=p(k)\lceil (n-k)/2\rceil-C_k(\varepsilon)$ with
$\varepsilon\equiv n-k\pmod 2$, valid for $n\ge4k+1$.}
\label{tab:PF}
\end{table}

\begin{remark}
All statements were also checked by computer. Sections~\ref{sec:struct}--\ref{sec:tree}, including the staircase
theorem and its consequences, were tested on all
$2^{16}-1$ numerical sets with Frobenius number at most $16$; Theorem~\ref{thm:twosizes} on all $(A^p,B^q)$ with
$A\le15$ and $p,q\le8$; and Table~\ref{tab:small} on all partitions with at most four parts and largest part at most
$12$. Theorems~\ref{thm:PG}, \ref{thm:PM}, \ref{thm:PMboundary}, \ref{thm:PMboundary2}, \ref{thm:PF} and
\ref{thm:PFboundary} were compared with a direct enumeration of all $51{,}861$ numerical semigroups whose partition has size at most $44$.
There were no exceptions.
\end{remark}

\section{Summary of solutions to the open problems of Burson, Nam and Sisneros-Thiry}\label{sec:solutions}

In this section we list the open problems of \cite{BNS} and give their solutions, all obtained from creation
sequences. Throughout, ``for all $n$ beyond a linear bound'' means for all $n$ at least a linear function of the
shift, and we state the bound explicitly each time.

\subsection{The common method}

Each problem asks about partitions close to an extreme case: many parts, many parts equal to $1$, or a hook that
contains almost all cells. In each case the creation sequence (equivalently, the threshold graph) of such a partition
consists of, an initial run of isolated vertices, the creation sequence of a small tail partition ($\nu$ or $\kappa$), a long run of dominating vertices, and a final isolated vertex.

By Lemma~\ref{lem:key} and Remark~\ref{rem:keygraph} (and by the corresponding step in the proof of
Theorem~\ref{thm:PG}, where the only critical vertex is $v_{2m}$), whether the partition corresponds to a semigroup
depends only on the creation sequence of the tail and on one integer offset. The offset is determined by the shift
($i$, $j$ or $k$) and not by $n$. Therefore the number of exceptions is a function of the shift alone, and it counts tails
with a prescribed arithmetic property of their creation sequences. This is exactly the kind of function that
\cite{BNS} asked for.

\subsection{Problem 1: genus (\cite[Remark 1]{BNS})}

\paragraph{Problem.} Note that $p(n)-\PG(3n-i,2n-i)$ and $b(n)-\PG(3n+2-i,2n-i)$ can be found for each
specific $i$, and state: ``it is an open problem to find functions $h_1(i)$ and $h_2(i)$ such that
$p(n)-\PG(3n-i,2n-i)=h_1(i)$ and $b(n)-\PG(3n+2-i,2n-i)=h_2(i)$ for all possible values of $i$.''

\paragraph{Solution.} By Theorem~\ref{thm:PG} and Corollary~\ref{cor:rem1},
\[
h_1(i)=E(i)\quad(n\ge i),\qquad h_2(i)=E(i+4)-4\quad(n\ge i+2),
\]
where $E(j)$ is the number of partitions $\nu$ such that $2|\nu|+\ell(\nu)-j$ is not a $1$-position of the creation
sequence of $\nu$ (Definition~\ref{def:E}). A closed generating function for $E$ is given in Proposition~\ref{prop:GF},
and its values are listed in Table~\ref{tab:EH}. In threshold-graph language, the tail $\nu$ is exceptional
exactly when the vertex $v_{2m}$ of $G(S_\lambda)$ is isolated.

\begin{example}
$h_1(4)=E(4)=4$. The tails counted are $\emptyset$ and $(1)$ (here $2|\nu|+\ell(\nu)<4$), together with $(1,1)$ and
$(3)$ (here $2|\nu|-\nu_t+t=4$ for $t=1$). For every $n\ge4$ the four partitions of $3n-4$ with $2n-4$ parts that do
not correspond to numerical semigroups are therefore
\[
(2^{n},1^{n-4}),\quad(3,2^{n-2},1^{n-3}),\quad(3,3,2^{n-4},1^{n-2}),\quad(5,2^{n-4},1^{n-1}),
\]
obtained from the tails $\emptyset,(1),(1,1),(3)$ by the parametrisation in the proof of Theorem~\ref{thm:PG}.
\end{example}

\subsection{Problem 2: multiplicity (\cite[Remark 4]{BNS})}

\paragraph{Problem.} ``It is an open problem to find a specific function $h(i)$ such that
$a(n+1)-\PM(2n-i,n-i)=h(i)$ for all possible values of $i$.''

\paragraph{Solution.} By Theorem~\ref{thm:PM} and Corollary~\ref{cor:rem4},
\[
h(i)=H(i)-1\qquad(n\ge2i+1),
\]
where $H(i)$ is the number of partitions $\kappa$ with $|\kappa|\le i$ such that $i-|\kappa|$ is a sum of two
$1$-positions of the creation sequence of $\kappa$ (Definition~\ref{def:H}). The term $-1$ removes the empty tail,
which corresponds to the one-part partitions already excluded from $a(n+1)$. At the boundary $n=2i$, where $3m=n$,
Theorems~\ref{thm:PMboundary} and \ref{thm:PMboundary2} give the exact corrections on the next two lines:
\[
a(2i+1)-\PM(3i,i)=H(i)-1-i,\qquad a(2i)-\PM(3i-2,i-1)=H(i)-1-(i-1).
\]
Below these lines only finitely many $n$ remain for each $i$.

\begin{example}
$h(3)=H(3)-1=5$. The tails with $|\kappa|\le3$ and $3-|\kappa|$ a sum of two $1$-positions are $(1)$ (offset $2$),
$(2)$ (offset $1$), and $(3)$, $(2,1)$, $(1,1,1)$ (offset $0$), plus the empty tail. The tail $(1,1)$ is not
exceptional: its creation sequence is $100$, and the offset $1$ is not a sum of two elements of $\{0,3,4,\dots\}$.
This corrects the range in \cite[Theorem 10-(4)]{BNS} to $n\ge7$ (Corollary~\ref{cor:PM}).
\end{example}

\subsection{Problem 3: Frobenius number (\cite[Remark 3]{BNS})}

\paragraph{Problem.} \cite[Remark 3]{BNS} give only the upper bound
$\PF(2n,2n-j)\le p(j)\lceil\frac{2n-j}{2}\rceil-\sum_{i=1}^{n-\lfloor j/2\rfloor}p(j,i)$, and note that the method
extends to further exceptional cases as $j$ grows.

\paragraph{Solution.} By Theorem~\ref{thm:PF} and Corollary~\ref{cor:rem3},
\[
\PF(2n,2n-j)=p(j)\Bigl\lceil n-\frac j2\Bigr\rceil-C_j(j\bmod2)\qquad(2n\ge4j+1),
\]
where $C_j(\varepsilon)=\sum_{\pi\vdash j}\bigl(c_\pi-G_\pi(\varepsilon)\bigr)$. Here $c_\pi=\pi_1+\ell(\pi)$ is the
length of the creation sequence of $\pi$, and $G_\pi(\varepsilon)$ counts the offsets $e<c_\pi$ of parity
$\varepsilon$ that are not sums of two $1$-positions of that sequence. More generally,
$\PF(n,n-k)=p(k)\lceil(n-k)/2\rceil-C_k(n-k)$ for $n\ge4k+1$. At the boundary $n=4k$ the exact correction is
$-\lceil k/2\rceil$ (Theorem~\ref{thm:PFboundary}), coming from the hooks of size $k$ whose creation sequences
$1\,0^{j}\,1^{k-j-1}\,0$ have too long a zero block (Lemma~\ref{lem:hooksums}).

\subsection{Problem 4: uniform theorems (\cite[Section 6]{BNS})}

\paragraph{Problem.} ``Many of the results above involve functions of the form $\PG(L_1(n)-k,L_2(n)-k)$,
$\PM(L_3(n)-k,L_4(n)-k)$, and $\PF(L_5(n),L_6(n)-k)$ for some specific linear functions $L_i(n)$ and certain values of
$k$. It is an open problem to extend the results uniformly to general theorems that hold for any value of
$k$.''

\paragraph{Solution.} The following three theorems hold for every value of the shift simultaneously.
\begin{align*}
\PG(n,g)&=p(n-g)-E(2n-3g) && \text{for } 2g\ge n && (\text{Theorem~\ref{thm:PG}}),\\
\PM(n,m)&=\bar p(n-m+1)-H(n-2m) && \text{for } m-1\le n<3m && (\text{Theorem~\ref{thm:PM}}),\\
\PF(n,n-k)&=p(k)\lceil (n-k)/2\rceil-C_k(n-k) && \text{for } n\ge4k+1 && (\text{Theorem~\ref{thm:PF}}).
\end{align*}
Along any line of the kind considered in \cite{BNS}, for instance
$\PG(L_1(n)-k,L_2(n)-k)$ with $2L_2(n)-k\ge L_1(n)$, substituting gives an explicit formula for every $k$
(Corollary~\ref{cor:rem1}). All counting theorems of \cite[Sections 3--5]{BNS} are special cases of these
theorems and of the boundary Theorems~\ref{thm:PMboundary}, \ref{thm:PMboundary2} and \ref{thm:PFboundary}
(Corollaries~\ref{cor:PG}, \ref{cor:PM} and \ref{cor:PF}). The three ranges are sharp
(Corollary~\ref{cor:PG}(1) and Theorems~\ref{thm:PMboundary} and \ref{thm:PFboundary}).

\section{Concluding remarks}
The creation-sequence point of view reduces the counting problems of \cite{BNS} to a single mechanism. Close to the
extreme cases, a partition fails to correspond to a numerical semigroup only because of a bounded tail $\nu$ or
$\kappa$, and whether it fails is decided by the Keith--Nath set of that tail. This is why the exceptional counts
are independent of $n$. The same words explain the symmetry $\ell(\lambda)\ge\lambda_1$ and the tree structure.
The following questions remain open.
\begin{enumerate}
\item Extend Theorems~\ref{thm:PG}, \ref{thm:PM} and \ref{thm:PF} below $g=n/2$, $m=n/3$ and $n=4k$. The boundary
theorems~\ref{thm:PMboundary} and \ref{thm:PFboundary} suggest that each further step adds a correction term counted
by hooks, or by partitions close to hooks.
\item For each fixed shift our formulas hold for all $n$ beyond a linear bound, and
the first boundary case is exact. The remaining finitely many small $n$ for each shift ($g<n/2$, $3m<n$, $n<4k$) are
not covered by a single formula (for the multiplicity, $3m<n-1$). Theorems~\ref{thm:PMboundary} and \ref{thm:PFboundary} suggest that each further
step introduces a correction counted by hooks or near-hooks.
\item The density conjecture of Constantin, Houston-Edwards and Kaplan, $S'(n)/P(n)\to0$, which \cite{BNS}
recall as motivation, is not resolved. We prove the bound $\limsup S'(n)/P(n)\le\frac12$
(Theorem~\ref{thm:density}).
\item Find the asymptotic growth of $E(j)$ and $H(i)$.
%\item Improve the density bound of Theorem~\ref{thm:density} towards the conjecture $S'(n)/P(n)\to0$ of \cite{CHK},
%for instance by combining $\ell(\lambda)\ge\lambda_1$ with the block-length bound of Corollary~\ref{cor:blocks}.
\end{enumerate}

{\bf Acknowledgement:} The first author would like to take this opportunity to gratefully acknowledge the American Mathematical Society and the International Mathematical Union for their generous support of his visit to Philadelphia, USA, to attend ICM 2026, and for providing excellent facilities and a stimulating academic environment that fostered fruitful research discussions and valuable academic interactions.\\

\textbf{Declaration of competing interests :} The authors declare no conflicts of interest.\\

\textbf{Data Availability:} Data sharing is not applicable to this article as no datasets were generated or analysed during the current study.

\end{document}